%% file: FSI_NonLinea.tex
\documentclass[12pt]{amsart}
\usepackage{graphicx}
\usepackage{caption}
\usepackage{subcaption}
\usepackage{color}
\usepackage{amsmath,amsfonts}
\usepackage{algorithm}
\usepackage{algorithmic}
\usepackage{booktabs}
\theoremstyle{definition}

\theoremstyle{remark}

\newcommand{\R}{\mathbb{R}}  % The real numbers.

\begin{document}

%%
%% The title of the paper goes here.  Edit to your title.
%%

\title[FSI Simulation with Hyperelastic Models]{Numerical Simulation of Fluid-Structure Interaction Problems with Hyperelastic Models: A Monolithic Approach}

%%
%% Now edit the following to give your name and address:
%% 
\author{Ulrich Langer}
\address{Johann Radon Institute for Computational and Applied Mathematics (RICAM), Austrian Academy of Sciences, Altenberger Strasse 69, A-4040 Linz, Austria}
\email{ulrich.langer@ricam.oeaw.ac.at}
\urladdr{http://www.ricam.oeaw.ac.at/people/u.langer/}% Delete if not wanted.

%%
%% If there is another author uncomment and edit the following.
%%

\author{Huidong Yang}
\address{Johann Radon Institute for Computational and Applied Mathematics (RICAM), Austrian Academy of Sciences, Altenberger Strasse 69, A-4040 Linz, Austria}
\email{huidong.yang@oeaw.ac.at}
\urladdr{http://people.ricam.oeaw.ac.at/h.yang/}

%%
%% If there are three of more authors they are added in the obvious
%% way. 
%%

%%%
%%% The following is for the abstract.  The abstract is optional and
%%% if not used just delete, or comment out, the following.
%%%

\begin{abstract}
  In this paper, we consider a monolithic approach
  to handle coupled fluid-structure interaction problems
  with different hyperelastic models in an all-at-once manner.
  We apply Newton’s method in the outer iteration
  dealing with nonlinearities of the coupled system. 
  We discuss preconditioned Krylov sub-space, algebraic multigrid and
  algebraic multilevel methods for solving the linearized algebraic equations.
  Finally, we compare the results of the monolithic approach
  with those of the corresponding partitioned approach
  that was studied in our previous work.
\end{abstract}

%%
%%  LaTeX will not make the title for the paper unless told to do so.
%%  This is done by uncommenting the following.
%%

\maketitle

%%
%% LaTeX can automatically make a table of contents.  This is done by
%% uncommenting the following:
%%

%\tableofcontents

%%
%% A Theorem is stated by
%%

%\begin{thm} The square of any real number is non-negative.
%\end{thm}

%%
%% Its proof is set off by
%% 

%\begin{proof}
%\end{proof}

%%
%% A new section is started as follows:
%%

%%%%%%%%%%%%%%%%%%%%%%%%%%%%%%%%%%%%%%%%%%%%%%%%%%%%%%%%%%%%%%%%%%%%%%
\section{Introduction}
%%%%%%%%%%%%%%%%%%%%%%%%%%%%%%%%%%%%%%%%%%%%%%%%%%%%%%%%%%%%%%%%%%%%%%
Parallel to the development of the partitioned approach 
for the fluid-structure interaction (FSI) 
simulation (see, e.g., \cite{DS:06, SB09:00, UK08:00, Habchi2013306, NME:NME1792}), 
the monolithic one also attracts 
many interests in the last decade; see, e.g., 
\cite{ABXCC10, EPFL11, Badia20084216, NME:NME3001, YB06:00, Heil20041, Razzaq20121156, CX10}. 
Compare to the flexibility of the partitioned approach, where existing 
fluid and structure sub-problem solvers can be directly 
reused or adapted in an iterative manner, 
the monolithic one behaves more stable and robust 
by dealing with the coupled nonlinear 
FSI system in an all-at-once manner. Formally speaking, 
we apply Newton's method (see \cite{PD05:00}) in an outer iteration 
dealing with nonlinearities originated 
from the domain movements, convection terms, material 
laws, transmission conditions and stabilization parameters (that may depend 
on the solution itself); as a price to pay, 
at each Newton iteration, a large linearized system 
is to be solved efficiently. 

In the monolithic approach, 
the linearization of the nonlinear coupled system turns out to be a nontrivial task and 
requires tedious work on both the analytical derivation and 
computer implementation. One difficulty 
considered in this work results from the hyperelastic nonlinear material law 
as for the thick-walled artery with the media and adventitia layer 
(see \cite{Holzapfel00:00, Holzapfel06:00}), 
for which the second and fourth order tensors of the energy 
functional with respect to the right Cauchy-Green tensor demand 
heavy amount of computational effort 
in each Newton iteration; see, e.g., \cite{GAH00:00, JB08:00} for an 
introduction on the basic tools used to derive these quantities under 
the Lagrangian framework and e.g., 
\cite{CA14} for the simulation of such arterial tissues. 
Thanks to our previous work in \cite{ULHY13}, 
the linearization for the hyperelastic models tackled in a partitioned FSI solver 
is reused in this work. 
Another difficulty stems from the fluid domain movement handled by 
the Arbitrary-Lagrangian-Eulerian (ALE) method, where the fluid domain 
displacement is introduced as an additional variable; 
see, e.g., \cite{TH:81, LF99:00, JD04:00}. To formalize 
the derivative of the fluid sub-problem with respect to the 
fluid domain displacement, the domain mapping (see, e.g., \cite{Wick11:00}) 
and shape derivative calculus (see, e.g., \cite{YB06:00, FENA2005127}) 
are two typical robust approaches 
mainly considered so far. In the domain mapping 
approach, the fluid sub-problem is mapped to the one 
on the reference (initial) fluid 
domain via the ALE mapping, 
that matches the Lagrangian structure domain on the 
interface for all the time. Therefore,  the FSI transmission conditions are defined 
on the unchanged interface between the fluid and structure reference domains. 
By transforming the fluid sub-problem 
from the current domain (ALE framework) to the reference 
domain (Lagrangian framework), 
the fluid domain deformation gradient tensor and its determinant arise, 
which leads to a formulation similar to the one 
under the Lagrangian framework as usually adopted in 
continuum mechanics. Thus, for the fluid sub-problem, 
we follow the same approach to compute the directional 
derivative with respect to the fluid domain displacement 
(see related techniques in, e.g., \cite{GAH00:00, JB08:00}) 
as we used for the hyperelastic equations in \cite{ULHY13}. 
In the second approach 
based on a shape derivative technique (see, e.g., \cite{JS92}), 
the derivative of the fluid sub-problem is then 
evaluated by computing the directional derivative with respect to 
the change of geometry (a small perturbation) on the current domain; 
see also this technique employed by 
the partitioned Newton's method in \cite{DS:06, Yang11:00}. 

In addition to the effort on the linearization of the coupled 
nonlinear system, the monolithic solver 
requires the properly designed preconditioners and solvers (as inner iteration) 
for the linearized coupled FSI system at each Newton iteration 
and may demand more effort. 
In \cite{Razzaq20121156}, the preconditioned 
Krylov subspace method (see, e.g., \cite{YS03:00}) 
and geometrical multigrid 
method (see, e.g., \cite{HB03}) with a Vanka-like smoother 
are employed to solve the linearized and discretized 2D FSI 
system using the high order $Q_2-P_1$ stabilized finite element pair. 
For the complex 3D geometries and unstructured meshes, 
in \cite{NME:NME3001}, the GMRES method (see \cite{Saad86}) 
accelerated by the block Gauss-Seidel preconditioner is considered, 
for which the block inverse is approximated 
by smoothed aggregation multigrid 
(see, e.g., \cite{SMTR08}) for each sub-problem. In order 
to improve the performance, 
a monolithic FSI algebraic multigrid (AMG) method using 
preconditioned Richardson iterations with potentially 
level-dependent damping parameters as smoothing 
steps is further developed therein. 
Besides, the monolithic solver is shown capable of utilizing 
parallel computing resources. In \cite{EPFL11}, parallel 
preconditioners of the coupled problem based on 
the algebraic additive Schwarz (see, e.g., \cite{ATOW05}) 
preconditioners for 
the sub-problems are built for both the convective explicit and 
geometry-convective explicit time discretized FSI systems. 
As a 2D counterpart, in \cite{ABXCC10}, a one-level additive 
Schwarz preconditioner (see, e.g., \cite{ATOW05}) for the linearized 
system is considered for the fully implicit time discretized FSI system, 
that is based on a sub-domain preconditioner constructed 
on an extention of a non-overlapping sub-domain to its neighbors. 

In this work, we focus on the 
development and comparison of different monolithic solution methods, 
namely, the Krylov 
subspace methods preconditioned by the block $LU$ decomposition of 
the coupled system, 
the AMG and algebraic multilevel 
(AMLI \cite{AO89I, AO90, PSV08:00, OA96, KJMS13, KJ12}, 
also referred to as K-cycle \cite{NY08, NLA542}) method, applied to the 
coupled FSI system with nearly incompressible 
hyperelastic models (see \cite{Holzapfel00:00,Holzapfel06:00}). 
Our solution methods are mainly based on 
a class of special AMG methods developed in 
\cite{FK98:00} and \cite{WM04:00, WM06:00}, 
for the discrete elliptic and saddle point problems, respectively, 
where the robust matrix-graph based coarsening strategies 
are proposed in a (pure) algebraic manner. This class of AMG methods 
have been applied to the 
sub-problems in the fluid-structure interaction simulation; 
see \cite{Yang11:00, YH11:00, HY10:000, ULHY13}. 
Particularly in our recent work \cite{ULHY13}, we 
have developed this approach by carefully 
choosing the effective smoothers: Braess-Sarazin 
smoother (see \cite{Braess97:00, WZ00:00}) and Vanka smoother 
(see \cite{Vanka86:00, TS09:00}), 
for the linearized Navier-Stokes equations under the ALE framework and 
hyperelastic equations under the Lagrangian 
framework, respectively. In order to further extend this class of AMG 
methods to the monolithic coupled FSI system after 
linearization, the two essential components in the AMG methods, 
the coarsening strategy and the smoother, for the 
coupled system are to be developed. Namely, 
the robust coarsening strategy using the stabilized 
Galerkin projection is constructured based on the $\inf-\sup$ 
condition (see, e.g., \cite{BF91:00, VG86}) on coarse levels 
for the indefinite sub-problems. 
By this means, we obtain the stabilized coupled systems on coarse levels. 
The effective smoother is designed by damped block Gauss-Seidel 
iterations applied to the coupled system,  that are based on the AMG cycles 
for the mesh movement, fluid and structure sub-problem, 
respectively. According to our numerical experiments, we observe 
the robustness of the damping parameter 
with respect to the AMG levels and different hyperelastic models adopted 
in the FSI simulation. As a variant of our coupled AMG method, 
we further consider the AMLI method for the coupled FSI system, in 
which we use the hierarchy of the coupled systems constructured in an algebraic 
manner as in the AMG methods. The smoothing for the 
coarse grid correction equation 
is performed by a flexible GMRES (FGMRES \cite{Saad1993}) scheme 
preconditioned by the multilevel preconditioner; see, e.g., \cite{PSV14} 
the application for the non-regularized 
Bingham fluid problem using the geometric multigrid method. In order to 
improve the performance, we finally consider the GMRES 
and FGMRES Kyrlov sub-space method 
preconditioned with such AMG and AMLI cycles. 

The remainder of the paper is organized in the following way. In 
Section \ref{sec:pre}, the coupled FSI system using a family of 
hyperelastic models for a model problem is formulated in a 
monolithic way. Section \ref{sec:dis} deals 
with the temporal and spatial discretization, and Newton's method 
tackling the linearization for the coupled nonlinear FSI system. 
In Section \ref{sec:lsm}, several monolithic solution 
methods for the linearized FSI system are considered in detail. 
Some numerical experiments are presented in Section \ref{sec:num}. 
Finally, in Section \ref{sec:con}, some conclusions are drawn. 

\section{A model problem}\label{sec:pre}
\subsection{Computational domains and mappings}
\begin{figure}[htbp]
  \centering
  \scalebox{0.70}{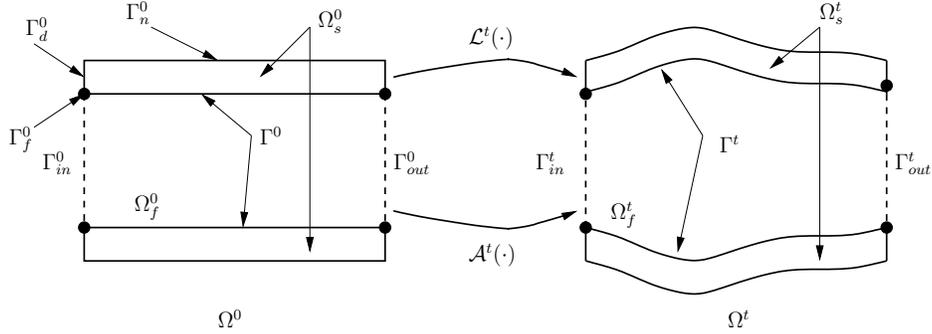}
  \caption{A schematic illustration of the domain mappings.}
  \label{fig:dom}
\end{figure}

We consider a model problem in the computational 
FSI domain $\Omega^t$ at time $t$ decomposed into the 
fluid domain $\Omega_f^t$ and 
the structure domain $\Omega_s^t$, i.e., 
$\overline{\Omega^t}=\overline{\Omega_f^t}\cup\overline{\Omega_s^t}$ and 
$\Omega_f^t\cap\Omega_s^t=\O$. Let 
$\Gamma_d^0$ and $\Gamma_n^0$ denote the boundaries 
with the homogeneous Dirichlet and Neumman condition for the 
structure sub-problem, respectively, 
$\Gamma_{in}^t$ and $\Gamma_{out}^t$ the boundaries with the inflow and outflow 
condition for the fluid sub-problem, 
respectively, $\Gamma_f^0=\overline{\Gamma_d^0}
\cap(\overline{\Gamma_{in}^0}\cup\overline{\Gamma_{out}^0})$ 
the fluid boundary with the homogeneous velocity condition, $\Gamma^t$ the 
interface between two domains: 
$\Gamma^t=\partial\Omega_f^t\cap\partial\Omega_s^t\setminus\Gamma_f^0$. 
At time $t=0$, we have 
all the initial configurations. See an illustration in Fig. \ref{fig:dom}. 

As usually adopted, we use the Lagrangian mapping 
(see, e.g., \cite{GAH00:00, JB08:00, WP08}) 
${\mathcal L}^t(\cdot):{\mathcal L}^t(x_0)=x_0+d_s(x_0, t)$ for 
all $x_0\in\Omega_s^0$ and $t\in(0, T)$ to track the motion 
of the structure body, where $d_s(\cdot, \cdot)$ denotes 
the structure displacement, i.e., 
$d_s(\cdot, \cdot):\Omega_s^0\times(0, T]\mapsto\R^3$. 
For the fluid sub-problem, we employ the 
arbitrary Lagrangian Eulerian (ALE, see, e.g., \cite{TH:81, JD04:00, LF99:00}) 
mapping ${\mathcal A}^t(\cdot)$ on $\Omega_f^0$ to 
track the fluid domain monition, i.e., 
${\mathcal A}^t(x_0)=x_0+d_f(x_0, t)$ for all 
$x_0\in\Omega_f^0$ and $t\in(0, T)$, where the 
fluid displacement $d_f(\cdot, \cdot):\Omega_f^0\times(0, T]\mapsto\R^3$ 
follows the fluid and structure particle motion on the $\Gamma^0$, 
and is arbitrary extended into the fluid domain $\Omega_f^0$ (see, e.g., 
\cite{Wick11:00}). The fluid domain velocity 
$w_f(\cdot, \cdot):\Omega_f^0\times(0, T]\mapsto\R^3$ is then given by 
$w_f=\partial_t d_f$. With help of this mapping, the fluid 
velocity $u(\cdot, \cdot):\Omega_f^0\times(0, T]\mapsto\R^3$ 
and pressure $p(\cdot, \cdot):\Omega_f^0\times(0, T]\mapsto\R$ 
are defined by the transformation:
\begin{subequations}\label{eq:vp}
  \begin{eqnarray}
    u(x, t)=\tilde{u}(\tilde{x}, t)=\tilde{u}({\mathcal A}_f^t(x), t),\\
    p(x, t)=\tilde{p}(\tilde{x}, t)=\tilde{p}({\mathcal A}_f^t(x), t),
  \end{eqnarray}
\end{subequations}
for all $x\in\Omega_f^0$ and $\tilde{x}={\mathcal A}_f^t(x)\in\Omega_f^t$. 
Here for simplicity 
of notations, $\tilde{u}(\cdot, \cdot)$ and $\tilde{p}(\cdot, \cdot)$ 
are used to indicate the variables under the Eulerian framework. 
\subsection{Basic notations in Kinematics}
In order to formulate the coupled system on the reference domain $\Omega^0$, 
we first introduce the following basic notations in Kinematics of nonlinear 
continuum mechanics; see, e.g., \cite{GAH00:00, JB08:00, WP08}. Let 
$F_f=\partial{\mathcal A}^t/{\partial x}=I+\nabla d_f$ for $x\in\Omega^0_f$ and 
$F_s=\partial{\mathcal L}^t/{\partial x}=I+\nabla d_s$ for $x\in\Omega^0_s$ 
denote the fluid and structure deformation gradient tensor, respectively. 
The determinant is given by $J_f=\text{det}F_f$ and $J_s=\text{det}F_s$, 
respectively. For the nonlinear hyperelastic models, further notations are used, 
namely, the right Cauchy-Green tensor $C=F_s^TF_s$ and the 
three principal invariants $I_1=C:I$, $I_2=0.5(I_1^2-C:C)$ and $I_3=\text{det}C$, 
respectively. Furthermore, the second Piola-Kirchoff tensor $S$ is defined by 
$S=2\partial\Psi/\partial C$, where $\Psi$ denotes the invariant dependent 
energy functional determined by the material properties. 
\subsection{A family of hyperelastic models}\label{sec:hypermat}
A family of hyperleastic models are used in this work, that posses nearly incompressible 
or anisotropic properties; see, e.g., 
\cite{JB08:00, GAH00:00, Holzapfel00:00, Holzapfel06:00, CA14}. 
We first consider the model of Neo-Hookean 
material, for which the energy functional is given by 
\begin{equation}\label{eq:nhenerg}
  \Psi=0.5c_{10}(J_1-3)+0.5\kappa(J_s-1)^2, 
\end{equation}
where $J_1=I_1I_3^{-1/3}$ denotes the invariant, $\kappa$ the bulk modulus and 
$c_{10}$ the material parameter related to the shear modulus. 
The second Piolad-Kirchoff tensor is then given by 
\begin{equation}\label{eq:nhpk}
  S=S^{'} -p_sJ_sC^{-1}
\end{equation}
with $S^{'}=c_{10}\partial J_1/\partial C$, where the structure pressure 
$p_s:=p_s(x, t)=-\kappa(J_s-1):\Omega_s^0\times(0, T]\mapsto\R^3$ is 
introduced in order to overcome the locking phenomena with large 
bulk modulus; see, e.g., 
\cite{Klaas96:00, Maniatty02:00, Goenezen11:00, TTE08:00}. 

We then consider the modified model of 
Mooney-Rivlin material, for which the energy functional is given by 
\begin{equation}\label{eq:mrenerg}
  \Psi=0.5c_{10}(J_1-3)+0.5c_{01}(J_2-3)+0.5\kappa(J_s-1)^2, 
\end{equation}
where $J_2=I_2I_3^{-2/3}$ denotes the invariant and $c_{01}$ the 
material parameter related to the shear modulus. 
The second Piolad-Kirchoff tensor is then accordingly given by 
\begin{equation}\label{eq:mrpk}
  S=S^{'}-p_sJ_sC^{-1},
\end{equation}
where $S^{'}=c_{10}\partial J_1/\partial C +c_{01}\partial J_2/\partial C$.

We finally consider the model of  the anisotropic 
two-layer thick-walled artery; 
see \cite{Holzapfel00:00, Holzapfel06:00}. The energy functional 
of such an arterial model is defined by 
\begin{equation}\label{eq:gaenerg}
  \Psi=0.5c_{10}(J_1-3)+\Psi_{\textup{aniso}}+0.5\kappa(J_s-1)^2
\end{equation}
with $\Psi_{\text{aniso}}=0.5k_1k_2\sum_{i=4, 6}(\exp(k_2(J_i-1)^2)-1)$, 
where $k_1$ and $k_2$ are a stress-like material parameter and 
a dimensionless parameter, respectively, 
associated with contribution of collagen to the response, and $J_4>1$ 
and $J_6>1$ are invariants active in extension, that are defined as 
$J_4=I_3^{-1/3}A_1:C$ and $J_6=I_3^{-1/3}A_2:C$, respectively. The 
tensors $A_1=a_{01}\otimes a_{01}$ and $A_2=a_{02}\otimes a_{02}$ are 
prescribed with the direction vectors 
$a_{01}=(0, \cos\alpha, \sin\alpha)^T$  and $a_{02}=(0,\cos\alpha, -\sin\alpha)^T$, 
respectively, where $\alpha\in\{\alpha_M, \alpha_A\}$ 
denotes the angle between the collagen fibers and 
the circumferential direction in the media and adventitia, respectively. 
The second Piola-Kirchoff tensor $S$ for this 
hyperelastic model is computed as
\begin{equation}\label{eq:gapk}
  S=S^{'}-p_sJ_sC^{-1}
\end{equation}
with $S^{'}=c_{10}\partial J_1/\partial C+k_1\sum_{i=4, 6}
(\exp(k_2(J_i-1)^2(J_i-1))\partial J_i/\partial C)$. In our numerical 
experiments, we use the geometrical configuration 
and material parameters of a rabbit carotid artery prescribed in 
\cite{Holzapfel00:00}. For modeling arterials in 
the FSI simulation considering specific fiber orientation, 
prestress and zero-stress configurations, 
and viscoelastic support conditions, we further refer to 
\cite{YB06:00, TET13:00, YB11:00, GJF12:00, GJF13:00, 
CNM12:00, GPlank12:00} for relevant details. 
\subsection{Coupled fluid-structure interaction system}
The coupled FSI system in strong form reads: Find $(d_f, u, p_f, d_s, p_s)$ 
such that 
 \begin{subequations}\label{eq:fsi}
    \begin{eqnarray}
      -\Delta d_f=0&{\textup{ in }} \Omega_f^0, \\
      d_f=d_s&{\textup{ on }} \Gamma^0,\\ [0.2cm]
      \rho_fJ_f\partial_t u+\rho_fJ_f((u- w_f)\cdot F_f^{-1}\nabla) u \\ \nonumber 
      -\nabla\cdot (J_f\sigma_f(u, p_f)F_f^{-T})=0&{\textup{ in }}\Omega_f^0,\\
      \nabla\cdot (J_fF_f^{-1}u)=0&{\textup{ in }}\Omega_f^0, \\[0.2cm]
      \rho_s\partial_{tt}d_s-\nabla\cdot (F_sS)=0&{\textup{ in }}\Omega_s^0,\\  
      -(J_s-1)-(1/\kappa)p_s=0&{\textup{ in }}\Omega_s^0, \\[0.2cm]
      u=\partial_t d_s&{\textup{ on }} \Gamma^0, \\
      J_f\sigma_{f}(u, p_f)F_f^{-T}n_f+F_sSn_s=0&{\textup{ on }}\Gamma^0,
    \end{eqnarray}
  \end{subequations}
 supplemented with the corresponding boundary conditions $d_f=0$ on 
$\Gamma^0_{in}\cup\Gamma^0_{out}$, $u=0$ on $\Gamma_f^0$, 
$J_f\sigma_f(u, p_f)F_f^{-T}n_f = g_{in}$ (a given function) on $\Gamma_{in}^t$ and 
$J_F\sigma_f(u, p_f)F_f^{-T}n_f = 0$ on $\Gamma_{out}^t$, $d_s=0$ on $\Gamma_d^0$ and 
$F_sSn_s=0$ on $\Gamma_n^0$, and proper initial conditions $u(x, 0)=0$ for all 
$x\in\Omega_f^0$ and $d_s(x,0)=\partial_td_s(x, 0)=0$ for all $x\in\Omega_s^0$. 
Here $\rho_f$ and $\rho_s$ denote the fluid and structure density, respectively, 
$n_f$ and $n_s$ the fluid and structure outerward unit normal vector, respectively, 
$\sigma_f(u, p_f):=\mu(\nabla u+\nabla^T)-p_fI$ the Cauchy stress tensor with the 
dynamic viscosity term $\mu$. Note that for the fluid sub-problem, we transform the 
momentum balance and mass conservation equations from the Eulerian to 
Lagrangian framework using the ALE mapping. 

\section{Temporal and spatial discretization and linearization}\label{sec:dis}
\subsection{The temporal discretization}
For the time discretization of the fluid sub-problem, we use 
the first order implicit Euler scheme. Let $u^n:=u(\cdot, t^n)$ 
and $w^n:=w(\cdot, t^n)=\partial_t d_f(\cdot, t^n)$ denote 
the approximateions of the fluid and fluid domain velocity at time level 
$t^n=n\Delta t$, $n=1,..., N$, $\Delta t=T/N$, i.e., the time period $(0, T]$ 
is subdivided into $N$ equidistant intervals. At time level $t^0$, 
the FSI solution is given by the initial conditions. The time derivatives at 
the level $t^n$ are then approximated as 
 \begin{subequations}\label{eq:tdfluid}
    \begin{eqnarray}
      \partial_t u^{n}&\approx&(u^{n}-u^{n-1})/\Delta t, \\
      w^{n}&\approx&(d^{n}-d^{n-1})/\Delta t.
    \end{eqnarray}
  \end{subequations}
For the structure sub-problem, a first order Newmark-$\beta$ scheme is used 
(see \cite{NM59:00}), i.e., 
 \begin{subequations}\label{eq:tdstruc}
   \begin{eqnarray}
     \ddot{d}_s^{n}&\approx&\frac{1}{\beta\Delta t^2}(d_s^n-d_s^{n-1})
     -\frac{1}{\beta\Delta t}\dot{d}_s^{n-1}-(\frac{0.5}{\beta}-1)\ddot{d}^{n-1}_s, \\
     \dot{d}_s^{n}&\approx&\dot{d}_s^{n-1}+\gamma\Delta t \ddot{d}_s^n
     +(1-\gamma)\Delta\ddot{d}_s^{n-1},
    \end{eqnarray}
  \end{subequations}
where $0<\beta\leq 1$ and $0\leq \gamma\leq 1$. 
\subsection{The time semi-discretized weak formulation}
In order to find the finite element FSI solution on 
proper function spaces, we first formulate the weak 
formulation for the coupled system (\ref{eq:fsi}). We 
introduce the notations $H^1(\Omega_f^0)$, $H^1(\Omega_s^0)$ 
and $L^2(\Omega_f^0)$ for the standard Sobolev and Lebesgue 
spaces (see, e.g., \cite{AF03:00}) on $\Omega_f^0$ and $\Omega_s^0$, 
respectively.  Let $V_m:=H^1(\Omega_f^0)^3$ be the fluid domain 
displacement space, $V_f:=H^1(\Omega_f^0)^3$ and 
$Q_f:=L^2(\Omega_f^0)$ be the fluid velocity and 
pressure space, respectively. The function spaces $V_s$ and $Q_s$ 
for the structure displacement and pressure shall be properly chosen regarding 
the nonlinearities; see, e.g., \cite{JMB:76, PGC88:00}. Incorporating boundary 
conditions, we further introduce the following function space notations: 
$V_{m,D}^n:=\{v\in V_m : v=d_s^n \text{ on } \Gamma^0\}$, 
$V_{m,0}:=\{v\in V_m : v=0 \text{ on } \Gamma_{in}^0\cup\Gamma_{out}^0\}$ 
for the mesh movement sub-problem, 
$V_{f,0}:=\{v\in V_f : v=0 \text{ on } \Gamma_f^0\}$ for the 
fluid sub-problem, $V_{s, D}^n:=\{v\in V_s : v=0 \text{ on } \Gamma_d^0 
\;|\; v=u^n\Delta t+v^{n-1} \text{ on } \Gamma^0\}$, 
$V_{s, 0}:=\{v\in V_s : v=0 \text{ on } \Gamma_d^0\cup\Gamma^0\}$ for 
the structure sub-problem. 

The weak formulation for the coupled system (\ref{eq:fsi}) reads: Find 
$(d^n_f$, $u_f^n$, $p_f^n$, $d_s^n$, $p_s^n)$
$\in$$(V_{m, D}$, $V_{f,0}$, $Q_f$, $V_{s, D}^n$, $Q_s)$ such that 
for all $(v_m$, $v_f$, $q_f$, $v_s$, $q_s)$ $\in$ 
$(V_{m,0}$, $V_{f, 0}$, $Q_f$, $V_{s, 0}$, $Q_s)$
\begin{subequations}\label{eq:weakfsi}
    \begin{eqnarray}
      (\nabla d_f^n, \nabla v_m)_{\Omega_f^0} = 0&, \\
      \rho_f(J_f(u^n-u^{n-1})/\Delta t, v_f)_{\Omega_f^0}\\ \nonumber 
      +\rho_f(J_f((u^n- (d^n_f-d^{n-1}_f)/\Delta t)\cdot F_f^{-1}\nabla) u^n, v_f)_{\Omega_f^0} \\\nonumber 
      -\langle g_{in}, v_f\rangle_{\Gamma_{in}^0}
      +(J_f\sigma_f(u^n, p_f^n)F_f^{-T}, \nabla v_f)_{\Omega_f^0}=0&,\\
      -(\nabla\cdot (J_fF_f^{-1}u^n), v_f)_{\Omega_f^0}=0&,\\
      (\beta_2d_s^n, v_s)_{\Omega_s^0}+(S^{'}, F_s^T\nabla v_s)_{\Omega_s^0}\\ \nonumber 
      -(J_s-1, q_s)_{\Omega_s^0}-(1/\kappa)(p_s^n, q_s)_{\Omega_s^0}=0&,\\ 
      -(p_sJ_sF_s^{-T}, \nabla v_s)_{\Omega_s^0}=0&, \\
      \langle J_f\sigma_{f}(u^n, p_f^n)F^{-T}n_f, v_f\rangle_{\Gamma^0}+
      \langle F_sSn_s, v_s\rangle_{\Gamma^0}=0&,
    \end{eqnarray}
  \end{subequations}
with $r_s=\beta_2d_s^{n-1}+\beta_1\dot{d}_s^{n-1}+\rho_s(0.5/\beta-1)\ddot{d}_s^{n-1}$, 
where $\beta_2=\rho_s/(\beta\Delta t^2)$ and $\beta_1=\rho_s/(\beta\Delta t)$. 
As observed, the fluid sub-problem is coupled with the mesh movement sub-problem 
in $\Omega_f^0$ and coupled with the structure sub-problem on $\Gamma^0$. 
The mesh movement sub-problem is coupled with the structure sub-problem on 
$\Gamma^0$. The equilibrium of surface tractions on $\Gamma^0$ is 
realized by the equilibrium of the 
residual of the weak formulation for the fluid and structure momentum equations with 
non-vanishing test functions $v_f\in V_{f,0}$ and $v_s\in V_{s,0}$ on $\Gamma^0$, 
where $v_f=v_s$ on $\Gamma^0$; see \cite{Yang11:00}.
\subsection{The spatial discretization and stabilization}
As in \cite{ULHY13}, we use Netgen \cite{JS97:00} to generate the tetrahedral mesh 
of the computational FSI domain $\Omega^0$ with conforming 
grids on the FSI interface and resolved different structure layers. 
For the mesh movement, we use $P_1$ finite element on the 
tetrahedral mesh. For the fluid and structure sub-problem, we use 
$P_1-P_1$ finite element with stabilization in order to fulfill the $\inf-\sup$ 
or LBB (Ladyshenkaya-Babu\v{s}ka-Brezzi) stability condition (see, e.g., 
\cite{BF91:00, VG86}), and to tackle the instability in 
advection dominated regions of the 
domain. In particular, we employ a unified 
streamline-upwind and pressure-stabilizing Petrov-Galerkin (SUPG/PSPG) 
method (see, e.g., \cite{Hughes89:00, Brooks82:00, WD06:00, FLD09:00, CF07:00}) 
to stabilize the $P_1-P_1$ discretized fluid sub-problem. 
For the structure sub-problem, we use the PSPG method 
(see, e.g., \cite{Hughes89:00, Klaas96:00, Maniatty02:00,Goenezen11:00}) 
to suppress the instability caused by equal order finite element interpolation 
spaces for the displacement and pressure. The application of this 
stabilization technique to the hyperelastic equations of anisotropic two-layer 
thick walled artery is prescribed in \cite{ULHY13}. 
\subsection{Newton's method for the nonlinear FSI system}\label{sec:nsm}
Formally speaking, after discretization in time and space 
of the coupled FSI system (\ref{eq:fsi}), 
we obtain the following nonlinear finite element algebraic equation 
\begin{equation}\label{eq:nlinsystem}
  {\mathcal R}(X) = 0
\end{equation}
with 
\begin{equation}
  {\mathcal R}(\cdot)=
  \left[\begin{array}{c}
      R_{ms}(\cdot)\\R_{mfs}(\cdot)\\R_{sf}(\cdot)
    \end{array}\right]\text{  and }
  X= 
  \left[\begin{array}{c}
      D_m \\ U_f\\ U_s
    \end{array}\right],
\end{equation}
where the subscripts $m$, $f$ and $s$ representing mesh movement, 
fluid and structure, respectively, and $ms$, $mfs$ and $sf$ the coupling 
among corresponding sub-problems. Furthermore 
$R_{ms}(X)=0$, $R_{mfs}(X)=0$ and $R_{sf}(X)=0$ stand for 
the finite element equations for the fluid mesh movement sub-problem coupled 
with the Dirichlet boundary condition on $\Gamma^0$ 
from the structure sub-problem, for the fluid sub-problem coupled 
with the fluid domain displacement from the mesh movement sub-problem 
and Neumann boundary condition on $\Gamma^0$ from the structure sub-problem, 
for the structure sub-problem coupled with the Dirichlet 
boundary condition on $\Gamma^0$ from the fluid sub-problem, respectively. The 
finite element solutions of the fluid domain displacement, fluid velocity and pressure, 
and structure displacement and pressure are denoted by $D_m$, $U_f$ 
and $U_s$, respectively. 

Newton's method applied to the nonlinear 
coupled FSI equation (\ref{eq:nlinsystem}) 
is presented in Algorithm \ref{alg:Newton}, where 
${\mathcal J}_k$ denote the Jacobian matrix and $\delta x_k$ the 
corrections of the finite element solutions at the $k$th step nonlinear iteration. 
\begin{algorithm}
   \caption{Newton's method}\label{alg:Newton} 
   Given initial $X_0$, for $k\ge 0$,
   \begin{algorithmic}[1]
     \STATE{assemble the Jacobian matrix ${\mathcal J}_k$ 
       at the current state $X_k$: for (\ref{eq:nlinsystem})
       \begin{equation}\label{eq:jaco}
         {\mathcal J}_k = {\mathcal R}^{'}(X_k)=
         \left[\begin{array}{ccc}
             \frac{\partial R_{ms}}{\partial D_m}& 
             0 &
             \frac{\partial R_{ms}}{\partial U_s}\\
             \frac{\partial R_{mfs}}{\partial D_m}& 
             \frac{\partial R_{mfs}}{\partial U_f}&
             \frac{\partial R_{mfs}}{\partial U_s}\\
             0&
             \frac{\partial R_{sf}}{\partial U_f}&
             \frac{\partial R_{sf}}{\partial U_s}\\
           \end{array}\right]
         \left(\begin{array}{c}
             D_{m,k}\\U_{f,k}\\U_{s,k}
           \end{array}\right),
       \end{equation}
     }
     \STATE{solve the linearize system 
       up to a relative residual error reduction factor $\varepsilon_l$:
       \begin{equation}\label{eq:jacosys}
         {\mathcal J}_k\delta x_k = -{\mathcal R}(X_k), 
       \end{equation}
     }
     \STATE{update the solution $X_{k+1}=X_k+\delta x_k$, 
       and go to step 1, until the relative residual error reduction 
       with the factor $\varepsilon_n$ is fulfilled:
       \begin{equation}\label{newstop}
         \|X_{k+1}-X_k\|_{L_2}
         \leq\varepsilon_n\|X_1-X_0\|_{L_2}.
       \end{equation}
     }
   \end{algorithmic}
\end{algorithm}
Note that the two terms 
$\partial R_{ms}/\partial U_s$ and 
$\partial R_{sf}/\partial U_f$ take 
the derivatives with respect to the 
structure displacement and fluid velocity, respectively, that 
corresponds to the linearization of 
two Dirichlet interface conditions on $\Gamma^0$ between 
the fluid and structure domain displacement, and between 
the fluid and structure velocity, respectively. The linearization 
of the Neumann interface condition on $\Gamma^0$ and of the 
fluid sub-problem are given 
in the second row of ${\mathcal J}_k$.  The linearization 
of structure sub-problem is given 
in the third row of ${\mathcal J}_k$. 

Besides the costly assembly procedure of ${\mathcal J}_k$ in 
(\ref{eq:jaco}), another main cost in Algorithm \ref{alg:Newton} is to 
solve the linearized equation (\ref{eq:jacosys}). 
More precisely, we come up with the linearized FSI system in the following 
reordered form (\ref{eq:linsystem1}) that we aim to solve 
(for simplicity of notations, we neglect the subscript $k$ and 
zero matrix entries): 
\begin{equation}\label{eq:linsystem1} 
  \left[\begin{array}{cccccccc}
      A_m^{ii}&A_m^{i\gamma}& & & & & & \\ [0.1cm] 
      & I & & & & -I& & \\ [0.1cm] 
      B^i_{fm}& B^\gamma_{fm}& -C_f & B^{i}_{1f} & B^{\gamma}_{1f} & & & \\ [0.1cm] 
      A^{ii}_{fm}&A^{i\gamma}_{fm}& B^{i}_{2f}&A^{ii}_f&A^{i \gamma}_f& & & \\ [0.1cm] 
      A^{\gamma i}_{fm}&A^{\gamma\gamma}_{fm}&B^{\gamma}_{2f}
      &A^{\gamma i}_f&A^{\gamma\gamma}_f
      &A^{\gamma\gamma}_s&A^{\gamma i}_s&B^{\gamma}_{2s}\\ [0.1cm] 
      & & & & -I & \frac{1}{\Delta t}I & & \\ [0.1cm] 
      & & & &    & A^{i\gamma}_s& A^{ii}_s& B^{i}_{2s}\\ [0.1cm] 
      & & & &    & B^{\gamma}_{1s}& B^{i}_{1s}& -C_s\\ 
    \end{array}
    \right]
  \left[\begin{array}{c}
      \Delta d_m^i\\ [0.1cm]  
      \Delta d_m^{\gamma} \\  [0.1cm] 
      \Delta p_f\\ [0.1cm]  
      \Delta u_f^i  \\   [0.1cm]  
      \Delta u_f^{\gamma}\\ [0.1cm] 
      \Delta d_s^{\gamma} \\  [0.1cm] 
      \Delta d_s^{i}\\ \Delta p_s
    \end{array}\right]   
  =
  \left[\begin{array}{c}
      r_m^i\\ [0.1cm] 
      r_m^{\gamma} \\ [0.1cm] 
      r_{p_f}\\ [0.1cm] 
      r_f^i \\[0.1cm] 
      r_f^{\gamma}\\[0.1cm]  
      r_s^{\gamma}\\[0.1cm] 
       r_s^{i}\\ [0.1cm] 
       r_{p_s}
    \end{array}\right], 
\end{equation}
where the superscripts $\gamma$ and $i$ are used to denote 
qualities associated to the nodal degrees of freedom (DOF) 
on the interface and the total remaining DOF in the domain 
and on the other boundaries of the domain. Furthermore, the 
qualities with the superscripts 
$\gamma\gamma$, $ii$, $\gamma i$ and $i\gamma$ indicate, 
that they result from the coupling of corresponding 
DOF. The solution posses a symbol $\Delta$ in front, 
indicating the DOF of the correction. It is easy to see from 
(\ref{eq:linsystem1}) how the sub-problems are linearized and 
coupled in a big FSI system. On the computer implementation, 
we are not explictly assembling the system (\ref{eq:linsystem1}), but the separate 
system for each sub-problem. The matching conditions are imposed implicitly 
by the conforming grids on the interface. 
This is easily implemented in the preconditioned Krylov subspace methods. 
The monolithic 
algebraic multigrid and multilevel approaches for the big coupled system require 
the formal systems on coarse levels. Furthermore, a direct solver is usually 
applied on the coarsest level, which requires an explictly formed system. 
Therefore, it is convenient and practical to form the big system in 
an explicit way and meanwhile to keep the flexibility of system assembling 
for each sub-problem. Therefore, we reformulate the system (\ref{eq:linsystem1}) as
\begin{equation}\label{eq:linsystem2}
  Kx=b
\end{equation}
\text{ with }
\begin{equation}\label{eq:linsystem3}
  K= 
  \left[\begin{array}{ccc}
      A_m& 0 & A_{ms}\\
      A_{fm}& A_{f}& A_{fs} \\
      0 &A_{sf}& A_s
    \end{array}
    \right],
  x=
  \left[\begin{array}{c}
      \Delta d_m
      \\ \Delta u_f
      \\ \Delta u_s
    \end{array}\right],
  b=  
  \left[\begin{array}{c}
      r_m
      \\ r_f
      \\ r_s
    \end{array}\right],
\end{equation}
where $A_m$, $A_f$ and $A_s$ represent 
the mesh movement, fluid and structure 
stiffness matrix from the finite element 
assembly, respectively , which are permuted from the 
corresponding ones in (\ref{eq:linsystem1}) according to their local 
nodal numbering of the finite element mesh for each sub-problem. The 
coupling matrix between $i\in\{m, f, s\}$ and $j\in\{m, f, s\}$ 
are denoted by $A_{ij}$, $i\neq j$. 
The solution vectors $\Delta d_m$, $\Delta u_f$ and $\Delta u_s$ denote 
the DOF of the correction for 
the fluid domain displacement, fluid velocity and pressure, 
and structure displacement and pressure, respectively. The residual 
vectors are presented by 
$r_m$, $r_f$ and $r_s$ for the fluid domain movement, fluid and structure sub-problem, 
respectively. These qualities are similar to the ones in (\ref{eq:linsystem1}),  
except that they are not recorded based on the 
separation of  the interface and remaing DOF. 
For consistency of notations, we will restrict our discussion 
to the solution methods of the linearized system (\ref{eq:linsystem2}) 
from now on. 

\section{Monolithic solution methods for the coupled system}\label{sec:lsm}
In this section, we discuss and compare different monolithic solution methods, 
namely, the preconditioned Krylov 
subspace methods, the algebraic multigrid and 
algebraic multilevel method, applied to the 
coupled system (\ref{eq:linsystem2}). 
\subsection{The preconditioned Krylov subspace methods}\label{sec:pky}
Because of the block structure of the system matrix $K$ in (\ref{eq:linsystem2}), 
we discuss some preconditioners mainly based on the $LU$ decomposition 
(see, e.g., \cite{YS03:00}). The inverse of 
the preconditioner applied to a given vector is easily realized using 
our efficient AMG methods for sub-problems (see \cite{ULHY13}). 
\subsubsection{The block-diagonal preconditioner} 
We first consider the simplest block-diagonal preconditioner $\tilde{P}_D$, 
that is obtained by approximating 
\begin{equation}\label{eq:pred}
  P_D = \left[\begin{array}{ccc}
      A_m & & \\
      & A_f & \\
      & & A_s \\
    \end{array}\right]
  \text{ with }
  \tilde{P}_D = \left[\begin{array}{ccc}
      \tilde{A}_m & & \\
      & \tilde{A}_f & \\
      & & \tilde{A}_s \\
    \end{array}\right],
\end{equation}
where $\tilde{A}_{i}=A_i(I-M_i^j)^{-1}$, $i\in\{m, f, s\}$, are corresponding 
multigrid preconditioners for each sub-problem; see, e.g, \cite{UH:02, UJ:89}. 
The inverse of $\tilde{P}_D$ is easily evaluated by 
\begin{equation}\label{eq:ipred}
  \tilde{P}_D^{-1} = \left[\begin{array}{ccc}
      \tilde{A}_m^{-1} & & \\
      & \tilde{A}_f ^{-1}& \\
      & & \tilde{A}_s^{-1} \\
    \end{array}\right],
\end{equation}
which corresponds to one AMG iteration applied to each sub-problem, that is 
developed in \cite{ULHY13}. 
This preconditioner completely neglects the coupling conditions among 
different sub-problems. 
\subsubsection{The block lower triangular preconditioner} 
The block lower triangular preconditioner $\tilde{P}_L$ 
is obtained by approximating 
\begin{equation}\label{eq:prel}
  P_L=\left[\begin{array}{ccc}
      A_m & &\\
      A_{fm} & A_f & \\
      0 & A_{sf}& A_s \\
    \end{array}\right]
  \text{ with }
  \tilde{P}_L=\left[\begin{array}{ccc}
      \tilde{A}_m & &\\
      A_{fm} & \tilde{A}_f & \\
       0 & A_{sf}& \tilde{A}_s \\
    \end{array}\right].
\end{equation}
It is easy to see the inverse of $\tilde{P}_L$ is given by 
\begin{equation}\label{eq:iprel}
  \tilde{P}_L^{-1}=\left[\begin{array}{ccc}
      \tilde{A}_m^{-1} & &\\
      -\tilde{A}_f^{-1}A_{fm}\tilde{A}_m^{-1} & \tilde{A}_f^{-1} & \\
      -\tilde{A}_{s}^{-1}A_{sf}\tilde{A}_f^{-1}A_{fm}\tilde{A}_m^{-1} & -\tilde{A}_s^{-1}A_{sf}\tilde{A}_f^{-1}& \tilde{A}_s^{-1} \\
    \end{array}\right],
\end{equation}
which is nothing but a block Gauss-Seidel 
iteration (using forward substitution) 
with zero initial guess applied to (\ref{eq:linsystem2}). 
This is easily computed since we have efficient AMG methods to 
approximate the inverse of $A_m$, $A_f$ and $A_s$. 
It is also easy to see one inverse operation 
of $\tilde{P}_L$ only requires (approximately) inverting $A_m$, 
$A_f$ and $A_s$ once. This 
preconditioner has taken into account the coupling block $A_{fm}$, 
the directional derivative of the fluid sub-problem 
with respect to the fluid domain displacement.
\subsubsection{The block upper triangular preconditioner} 
We then consider the block upper triangular preconditioner $\tilde{P}_U$ 
obtained by approximating  
\begin{equation}\label{eq:prer}
  P_U=\left[\begin{array}{ccc}
      A_m & 0 & A_{ms}\\
      & A_f &  A_{fs}\\
      & & A_s \\
    \end{array}\right]
  \text{ with }
  \tilde{P}_U=\left[\begin{array}{ccc}
      \tilde{A}_m & 0 & A_{ms}\\
      & \tilde{A}_f &  A_{fs}\\
      & & \tilde{A}_s \\
    \end{array}\right].
\end{equation}
The coupling blocks $A_{ms}$ and $A_{fs}$ are included, which 
represent the coupling of the Dirichlet interface condition between 
the fluid and structure domain displacement, and the coupling of the 
Neumann interface condition between the fluid and the structure 
surface traction, respectively. The inverse $\tilde{P}_U^{-1}$ is given by
\begin{equation}\label{eq:iprer}
  \tilde{P}_U^{-1}=\left[\begin{array}{ccc}
      \tilde{A}_m^{-1} & 0 & -\tilde{A}_{m}^{-1}A_{ms}\tilde{A}_s^{-1}\\
      & \tilde{A}_f^{-1} &  -\tilde{A}_f^{-1}A_{fs}\tilde{A}_s^{-1}\\
      & & \tilde{A}_s^{-1} \\
    \end{array}\right],
\end{equation}
that is nothing but a Gauss-Seidel iteration using a backward substitution. As we 
see the block $A_{fm}$ of the derivative of the fluid sub-problem 
with respect to the fluid domain displacement is not take into account. 
\subsubsection{The $SSOR-$preconditioner}
We consider a Symmetric Successive Over-Relaxation (SSOR) with a 
special choice of the relaxation parameter $\omega=1$. The preconditioner 
is based on the following $LU$ factorization of $P_{SSOR}$ given by 
\begin{equation}\label{eq:pres}
  \begin{aligned}
    P_{SSOR}&=\left[\begin{array}{ccc}
        A_m &  & \\
        A_{fm} &  A_f  & \\
        0 & A_{sf} & A_s \\
      \end{array}\right]\times
    \left[\begin{array}{ccc}
        I& 0& A_m^{-1}A_{ms}\\
        & I & A_f^{-1}A_{fs}\\
        & & I \\
      \end{array}\right]\\
    &=\left[\begin{array}{ccc}
        A_m &  0 & A_{ms}\\
        A_{fm} &  A_f  & A_{fs}+A_{fm}A_m^{-1}A_{ms} \\
        0 & A_{sf} & A_s+A_{sf}A_f^{-1}A_{fs} \\
      \end{array}\right].
  \end{aligned}
\end{equation}
that can be reformulated as $P_{SSOR}=K+R_{SSOR}$, where the remainder 
$R_{SSOR}$ is given by
$$
R_{SSOR}=\left[\begin{array}{ccc}
    0&  0 & 0\\
    0&  0 & A_{fm}A_m^{-1}A_{ms} \\
    0 & 0 & A_{sf}A_f^{-1}A_{fs} \\
  \end{array}\right].
$$
The $SSOR$ preconditioner $\tilde{P}_{SSOR}$ is then given by 
\begin{equation}
  \begin{aligned}
    \tilde{P}_{SSOR}&=\left[\begin{array}{ccc}
        \tilde{A}_m &  & \\
        A_{fm} &  \tilde{A}_f  & \\
        0 & A_{sf} & \tilde{A}_s \\
      \end{array}\right]\times
    \left[\begin{array}{ccc}
        I& 0& \tilde{A}_m^{-1}A_{ms}\\
        & I & \tilde{A}_f^{-1}A_{fs}\\
        & & I \\
      \end{array}\right].
  \end{aligned}
\end{equation}
The inverse of $\tilde{P}_{SSOR}$ is computed by two block 
Gauss-Seidel iterations using the backward and forward 
substitution consecutively: 
\begin{equation}\label{eq:ipres}
  \begin{aligned}
    \tilde{P}_{SSOR}^{-1}=
    &\left[\begin{array}{ccc}
        I & 0& -\tilde{A}_m^{-1}A_{ms}\\
         & I & -\tilde{A}_f^{-1}A_{fs}\\
        & & I \\
      \end{array}\right]\times
    \\
    &\left[\begin{array}{ccc}
        \tilde{A}_m^{-1} & &\\
        -\tilde{A}_f^{-1}A_{fm}\tilde{A}_m^{-1} & \tilde{A}_f^{-1} & \\
        -\tilde{A}_{s}^{-1}A_{sf}\tilde{A}_f^{-1}A_{fm}\tilde{A}_m^{-1} & -\tilde{A}_s^{-1}A_{sf}\tilde{A}_f^{-1}& \tilde{A}_s^{-1} \\
      \end{array}\right].
  \end{aligned}
\end{equation}
Compared to $\tilde{P}_L^{-1}$ and $\tilde{P}_U^{-1}$, two more inverse operations of $\tilde{A}_m^{-1}$ 
and $\tilde{A}_f^{-1}$ are required. 
\subsubsection{The $ILU(0)-$preconditioner}
We finally consider the $ILU(0)-$ preconditioner $\tilde{P}_{ILU}$. 
This incomplete factorization technique 
is described in, e.g., \cite{YS03:00, VH97:00, OA96}. Here we apply a 
block $ILU(0)$ factorization for the coupled FSI system given by
\begin{equation}\label{eq:preilu}
  \begin{aligned}
    P_{ILU}&=\left[\begin{array}{ccc}
        I &  & \\
        A_{fm}A_m^{-1} &  I  & \\
        0 & A_{sf}A_f^{-1} & I \\
      \end{array}\right]\times
    \left[\begin{array}{ccc}
        A_m& 0& A_{ms}\\
        & A_f & A_{fs}-A_{fm}A_m^{-1}A_{ms}\\
        & & A_s \\
      \end{array}\right]\\
    &=\left[\begin{array}{ccc}
        A_m &  0 & A_{ms}\\
        A_{fm} &  A_f  & A_{fs} \\
        0 & A_{sf} & A_s+A_{sf}A_f^{-1}(A_{fs}-A_{fm}A_m^{-1}A_{ms}) \\
      \end{array}\right],
  \end{aligned}
\end{equation}
that can be rewritten as $P_{ILU}=K+R_{ILU}$, where the 
remainder $R_{ILU}$ is given by 
$$
R_{ILU}=\left[\begin{array}{ccc}
    0&  0 & 0\\
    0&  0 & 0 \\
    0 & 0 & A_{sf}A_f^{-1}(A_{fs}-A_{fm}A_m^{-1}A_{ms}) \\
  \end{array}\right].
$$
The preconditioner $\tilde{P}_{ILU}$ is then given by 
\begin{equation}
  \begin{aligned}
    \tilde{P}_{ILU}&=\left[\begin{array}{ccc}
        I &  & \\
        A_{fm}\tilde{A}_m^{-1} &  I  & \\
        0 & A_{sf}\tilde{A}_f^{-1} & I \\
      \end{array}\right]\times
    \left[\begin{array}{ccc}
        \tilde{A}_m& 0& A_{ms}\\
        & \tilde{A}_f & A_{fs}-A_{fm}\tilde{A}_m^{-1}A_{ms}\\
        & & \tilde{A}_s \\
      \end{array}\right]\\
    &=\left[\begin{array}{ccc}
        \tilde{A}_m &  0 & A_{ms}\\
        A_{fm} &  \tilde{A}_f  & A_{fs} \\
        0 & A_{sf} & \tilde{A}_s+A_{sf}\tilde{A}_f^{-1}(A_{fs}-A_{fm}\tilde{A}_m^{-1}A_{ms}) \\
      \end{array}\right].
  \end{aligned}
\end{equation}
The inverse of $\tilde{P}_{ILU}$ is then computed by two block 
Gauss-Seidel iterations using the forward and backward 
substitution consecutively: 
\begin{equation}\label{eq:ipreilu}
  \begin{aligned}
    \tilde{P}_{ILU}^{-1}=
    &\left[\begin{array}{ccc}
        I & & \\
        -A_{fm}\tilde{A}_m^{-1}& I & \\
        A_{sf}\tilde{A}_f^{-1}A_{fm}\tilde{A}_m^{-1}& -A_{sf}\tilde{A}_f^{-1}& I \\
      \end{array}\right]\times
    \\
    &\left[\begin{array}{ccc}
        \tilde{A}_m^{-1} & 0 & -\tilde{A}_m^{-1}\tilde{A}_s^{-1}\\
        & \tilde{A}_f^{-1} &-\tilde{A}_f^{-1}(A_{fs}-A_{fm}\tilde{A}_m^{-1}A_{ms}) \\
        & &  \tilde{A}_s^{-1}\\
      \end{array}\right].
  \end{aligned}
\end{equation}
Compared to $\tilde{P}_L^{-1}$ and $\tilde{P}_U^{-1}$, 
two more inverse operations of $\tilde{A}_m^{-1}$ 
and one more inverse operation of $\tilde{A}_f^{-1}$ are required. 

\subsection{Algebraic multigrid method for the coupled FSI system}\label{sec:amgfsi}
The $LU$ factorization is probably the best well-known preconditioner for solving 
general systems. Unfortunately, those preconditioners discussed in Section 
\ref{sec:pky} for the FSI coupled system are not robust with respect to, e.g., 
the mesh size. As we observe from numerical experiments, the iteration numbers 
increase when the mesh is refined. In order to eliminate the mesh dependence, we 
consider the AMG and AMLI method. These methods 
tackles the low and high frequency errors separately by using the smoothing and 
coarse grid correction step. We discuss two essential components, the 
coarsening strategy and smoother, that are used in both the AMG and AMLI method.
\subsubsection{The coarsening strategy}
First of all, we define a full rank prolongation matrix 
\begin{equation}\label{eq:pro}
  P_{l+1}^l =\left[
  \begin{array}{ccc}
    P_m^l& &\\
    & P_f^l & \\
    & & P_s^l
  \end{array}\right],
\end{equation}
where $l=1,...,L-1$, indicates the levels of a hierarchy, i.e., index $1$ 
refers to the finest level and $L$ the coarsest level. Here 
$P_m^l:\R^{n_m^{l+1}}\mapsto\R^{n_m^l}$ denotes the prolongation matrix constructured 
for the elliptic mesh movement sub-problem as in \cite{FK98:00},  $n_m^l$ 
the number of DOF of the mesh movement sub-problem 
on level $l$ and $n_m^{l+1}<n_m^l$. In a similar way, 
$P_f : \R^{n_m^{l+1}}\mapsto\R^{n_m^l}$ and 
$P_s : \R^{n_s^{l+1}}\mapsto\R^{n_s^l}$ represent the prolongation matrices constructured 
for the indefinite fluid and structure sub-problem as in \cite{WM04:00, ULHY13}, 
that take the stability into account by proper scaling and 
avoid a mixture of velocity/displacement and 
pressure components on coarse levels, $n_f^l$ and $n_s^l$ 
the number of DOF of the fluid and structure sub-problem 
on level $l$ and $n_f^{l+1}<n_f^l$, $n_s^{l+1}<n_s^l$. Then it is easy to see 
$P_{l+1}^l : \R^{n_m^{l+1}+n_f^{l+1}+n_s^{l+1}}\mapsto\R^{n_m^l+n_f^l+n_s^l}$. More sophisticated 
and expensive coarsening strategies of the AMG method for 
saddle point systems arising from the fluid sub-problem can be found in \cite{BM:13}. 
In this work, we restrict ourselves to the strategy introduced in \cite{WM04:00}, 
where a simple scaling technique is applied. We then define a restriction matrix 
$R_l^{l+1} : \R^{n_m^l+n_f^l+n_s^l}\mapsto\R^{n_m^{l+1}+n_f^{l+1}+n_s^{l+1}}$ as 
\begin{equation}\label{eq:res}
  R_l^{l+1} =\left[
  \begin{array}{ccc}
    R_m^{l+1}& &\\
    & R_f^{l+1} & \\
    & & R_s^{l+1}
  \end{array}\right],
\end{equation}
where $R_m^{l+1}=(P_m^l)^T$, $R_f^{l+1}=(P_f^l)^T$ and $R_s^{l+1}=(P_s^l)^T$. 
The system on the finest level $L$ is given by (\ref{eq:linsystem2}) that 
is formulated as $K_1x_1=b_1$. Then the system matrix on the coarse level $l+1$ is 
formulated by the Galerkin projection that has considered the stability of indefinite 
sub-systems on coarse levels:
\begin{equation}\label{eq:cosys}
  K_{l+1}=R_l^{l+1}K_lP_{l+1}^l=\left[
  \begin{array}{ccc}
    R_m^{l+1}A_m^lP_m^l&0&R_m^{l+1}A_{ms}^lP_s^l\\
    R_f^{l+1}A_{fm}^lP_m^l&R_f^{l+1}A_f^lP_f^l& R_f^{l+1}A_{fs}^lP_s^l\\
    0&R_s^{l+1}A_{sf}^lP_f^l&R_s^{l+1}A_s^lP_s^l
  \end{array}\right],
\end{equation}
where $A^l_{i}$, $i\in\{m, ms, fm, f, fs, sf, s\}$ denote the matrices on the level $l$, 
$l=1,...,L-1$. On the coarsest level $L$, the coupled system is solved by a direct solver. 
\subsubsection{The smoother}\label{sec:fsism}
To complete the AMG method we need an iterative method (the smoother) 
for the problem $K_lx_l=b_l$, $l=1,...,L-1$,
\begin{equation}\label{eq:sm}
  x_l^{k+1}={\mathcal S}_l(x_l^k, b_l)
\end{equation}
with
\begin{equation}\label{eq:smvec}
  x_l^k=\left[
  \begin{array}{c}
    \Delta d_{m, l}^k\\ \Delta u_{f, l}^k\\ \Delta u_{s, l}^k
  \end{array}\right],
  b_l=\left[
    \begin{array}{c}
      r_{m, l}\\ r_{f, l}\\ r_{s, l}
    \end{array}\right],
\end{equation}
where $k$ is the iteration index.

For this coupled FSI system, we consider the following 
preconditioned Richardson method, that turns out to be 
an effective FSI smoother with sufficient large number of smoothing 
steps: For $k\geq 0$, 
\begin{equation}\label{eq:smrd}
  \left[
  \begin{array}{c}
    \Delta d_{m, l}^{k+1}\\
    \Delta u_{f, l}^{k+1}\\
    \Delta u_{s. l}^{k+1}
  \end{array}
  \right]
  =
  \left[
  \begin{array}{c}
    \Delta d_{m, l}^k\\
    \Delta u_{f, l}^k\\
    \Delta u_{s. l}^k
  \end{array}
  \right]
  +
  P_{Rich}^{-1}
  \left(
  \left[
  \begin{array}{c}
    r_{m, l}\\
    r_{f, l}\\
    r_{s. l}
  \end{array}
  \right]-
  K_l
  \left[
  \begin{array}{c}
    \Delta d_{m, l}^k\\
    \Delta u_{f, l}^k\\
    \Delta u_{s. l}^k
  \end{array}
  \right]
  \right),
\end{equation}
where the preconditioner is given by 
\begin{equation}\label{eq:prich}
  P_{Rich}=  
  \left[
    \begin{array}{ccc}
      \frac{1}{\omega_m}\tilde{A_m^l}&  & \\
      A_{fm}^l& \frac{1}{\omega_f}\tilde{A_f^l}& \\
      0& A_{sf}^l & \frac{1}{\omega_s}\tilde{A_s^l}\\
    \end{array}
    \right]
\end{equation}
with the scaled block diagonal matrices. The inverse of each of these 
matrices is realized by applying one AMG cycle to each sub-problem, that 
has been developed in our previous work \cite{ULHY13}. In principle, 
these damping parameters $\omega_i$, $i\in\{m, f, s\}$, 
may be chosen differently. For simplicity, 
we use $\omega_m=\omega_f=\omega_s=\omega$ in our numerical experiments. 
This FSI smoother shows numerical robustness with respect to different hyperelastic models 
considered in the coupled FSI system and the AMG levels, i.e., the same damping 
parameter $\omega$ has been used in our numerical experiments. 
It is easy to see, one iteration of the 
preconditiond Richardson method consists of three steps of 
the following damped block Gauss-Seidel like iteration, 
that is demonstrated in Algorithm \ref{alg:gs}. 
\begin{algorithm}
   \caption{Block Gauss-Seidel iteration: 
     $x_l^{k+1}={\mathcal S}_l(x_l^k, b_l)$ }\label{alg:gs} 
   Given initial $x_l^k$,
   \begin{algorithmic}[1]
     \STATE{$\Delta d_{m, l}^{k+1}=\Delta d_{m, l}^k
       +\omega_m (\tilde{A_m^l})^{-1}(r_{m, l}-A_m^l\Delta d_{m, l}^k
       -A_{ms}^l\Delta u_{s, l}^k)$,
     }
     \STATE{$\Delta u_{f, l}^{k+1}=\Delta u_{f, l}^k
       +\omega_f(\tilde{A_f^l})^{-1}(r_{f, l}-A_{fm}^l\Delta d_{m, l}^{k+1}
       -A_f^l\Delta u_{f, l}^k-A_{f, s}^l\Delta u_{s, l}^k)$,
     }
     \STATE{$\Delta u_{s. l}^{k+1}=\Delta u_{s, l}^k
       +\omega_s(\tilde{A_s^l})^{-1}
       (r_{s, l}-A_{sf}^l\Delta u_{f, l}^{k+1}-A_s^l\Delta u_{s, l}^k)$.
     }
   \end{algorithmic}
\end{algorithm}
\subsubsection{The algebraic multigrid iteration}
The basic AMG iteration is given in Algorithm \ref{alg:amgit}, 
where $m_{pre}$ and $m_{post}$ refer to the 
number of pre- and post-smoothing steps. For $\nu=1$ and 
$\nu=2$, the iterations in Algorithm \ref{alg:amgit} are called 
V- and W-cycle, respectively. In our numerical experiments, we choose 
the W-cycle. On the coarsest level $L$, we use direct solver to handle the 
coupled system. 
\begin{algorithm}
   \caption{Basic AMG iteration: AMG($K_l, x_l, b_l$) }\label{alg:amgit} 
   \begin{algorithmic}[1]
     \FOR{$k=1$ to $m_{pre}$} 
     \STATE $x_l^{k+1}={\mathcal S}_l(x_l^k, b_l)$
     \ENDFOR
     \STATE $b_{l+1}=R_l^{l+1}(b_l-K_lx_l)$,
     \IF{l+1=L}
     \STATE Solve $K_Lx_L=b_L$
     \ELSE
     \STATE{
       $x_{l+1}=0$,
       \FOR{$k=1,...,\nu$}
       \STATE $x_{l+1}=$AMG($K_{l+1}, x_{l+1}, b_{l+1}$),
       \ENDFOR
     }
     \ENDIF
     \STATE $x_l=x_l+P_{l+1}^lx_{l+1}$,
     \FOR{$k=1$ to $m_{post}$} 
     \STATE $x_l^{k+1}={\mathcal S}_l(x_l^k, b_l)$,
     \ENDFOR
     \STATE return $x_l$.
   \end{algorithmic}
\end{algorithm}

As seen from Algorithm \ref{alg:amgit}, steps 1-3 and steps 14-16 correspond to 
the presmoothing and postsmoothing, respectively, steps 4-13 are referred to 
as "coarse grid correction". The full AMG iterations are realized by repeated application 
of this algorithm until it satisfies certain stopping criteria. The iteration in this algorithm 
is also combined with GMRES \cite{Saad86} and FGMRES \cite{Saad1993} 
methods, that leads to fast convergence of the preconditioned Krylov subspace methods 
for the coupled FSI system. 
\subsection{Algebraic multilevel method for the coupled FSI system}\label{sec:amlgsi}
The AMLI method \cite{AO89I, AO90, PSV08:00, KJMS13}, 
sometimes referred to as "K-cycle", is viewed 
as a W-cycle with the Krylov acceleration at the 
intermediate levels; see, e.g. \cite{NY08, NLA542, PSV14, PSV08:00}. 
Here we combine our monolithic 
AMG method with the FGMRES Krylov subspace method 
at the intermediate levels, i.e., 
we reuse the coarsening strategy and 
smoothers constructed for the FSI AMG method. 
Instead of 
calling the AMG cycle (steps 9-10 in Algorithm \ref{alg:amgit}), the AMLI algorithm 
calls the AMLI cycle recursively $\nu$ times as a preconditioner inside the 
FGMRES method for the coarse grid correction equations; 
see step 9 in Algorithm \ref{alg:amliit}.
\begin{algorithm}
   \caption{Basic AMLI iteration: AMLI($K_l, x_l, b_l$) }\label{alg:amliit} 
   \begin{algorithmic}[1]
     \FOR{$k=1$ to $m_{pre}$} 
     \STATE $x_l^{k+1}={\mathcal S}_l(x_l^k, b_l)$
     \ENDFOR
     \STATE $b_{l+1}=R_l^{l+1}(b_l-K_lx_l)$,
     \IF{l+1=L}
     \STATE Solve $K_Lx_L=b_L$
     \ELSE
     \STATE{$x_{l+1}=0$,}
     \STATE{FGMRES($K_{l+1}$, $x_{l+1}$, $b_{l+1}$, $\nu$, AMLI),}
     \ENDIF
     \STATE $x_l=x_l+P_{l+1}^lx_{l+1}$,
     \FOR{$k=1$ to $m_{post}$} 
     \STATE $x_l^{k+1}={\mathcal S}_l(x_l^k, b_l)$,
     \ENDFOR
     \STATE return $x_l$.
   \end{algorithmic}
\end{algorithm}

\begin{figure}[htbp]
  \centering{
    \begin{subfigure}[b]{0.8\textwidth}
      \includegraphics[scale=0.4]{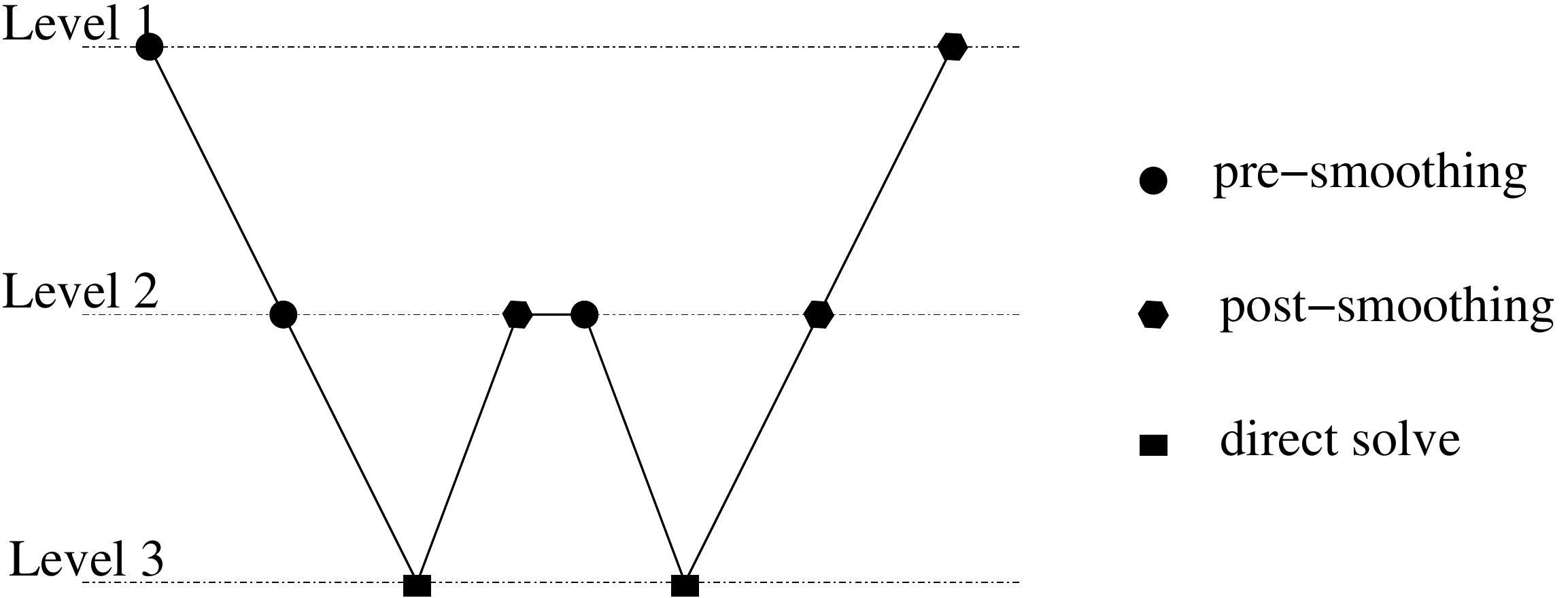}
      \caption{An AMG W-cycle}\label{fig:amli}
    \end{subfigure}
    \vfill
    \begin{subfigure}[b]{0.8\textwidth}
      \includegraphics[scale=0.4]{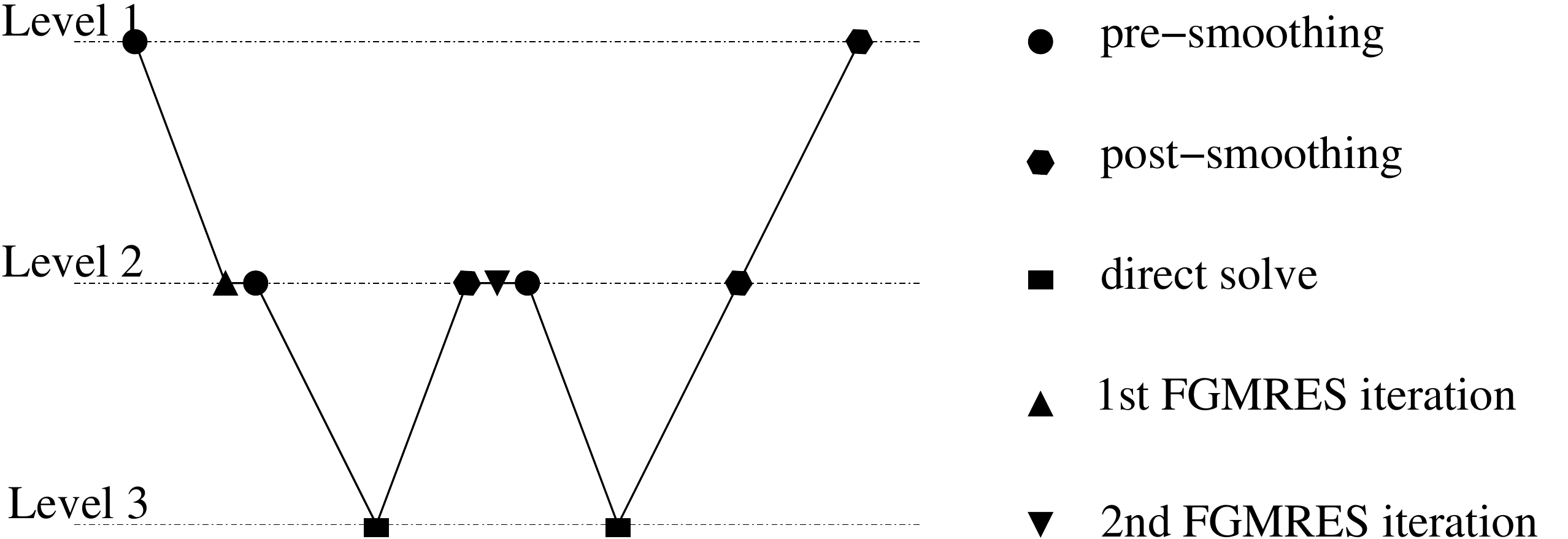}
      \caption{An AMLI W-cycle}\label{fig:amli}
    \end{subfigure}
  }
  \caption{An illustration of W-cycles for the AMG and AMLI methods with 3 levels.}
  \label{fig:amgamli}
\end{figure}

It is easy to see, this method represents a variant of the W-cycle AMG method in 
the case of $\nu=2$; see an illustration for such W-cyles 
with $3$ levels ($L=3$) in Fig. \ref{fig:amgamli}.  Compared to the AMG W-cycle, 
the two preconditioned FGMRES iterations are called consecutively on the second 
level of the AMLI W-cycle, that are used to accelerate the convergence rate.  

\section{Numerical experiments}\label{sec:num}
\subsection{Material and geometrical data, meshes and boundary conditions}
\begin{figure}[ht!]
  \centering{
    %\resizebox{5cm}{5cm}{\input{gageo.pdf_t}}
    \scalebox{0.30}{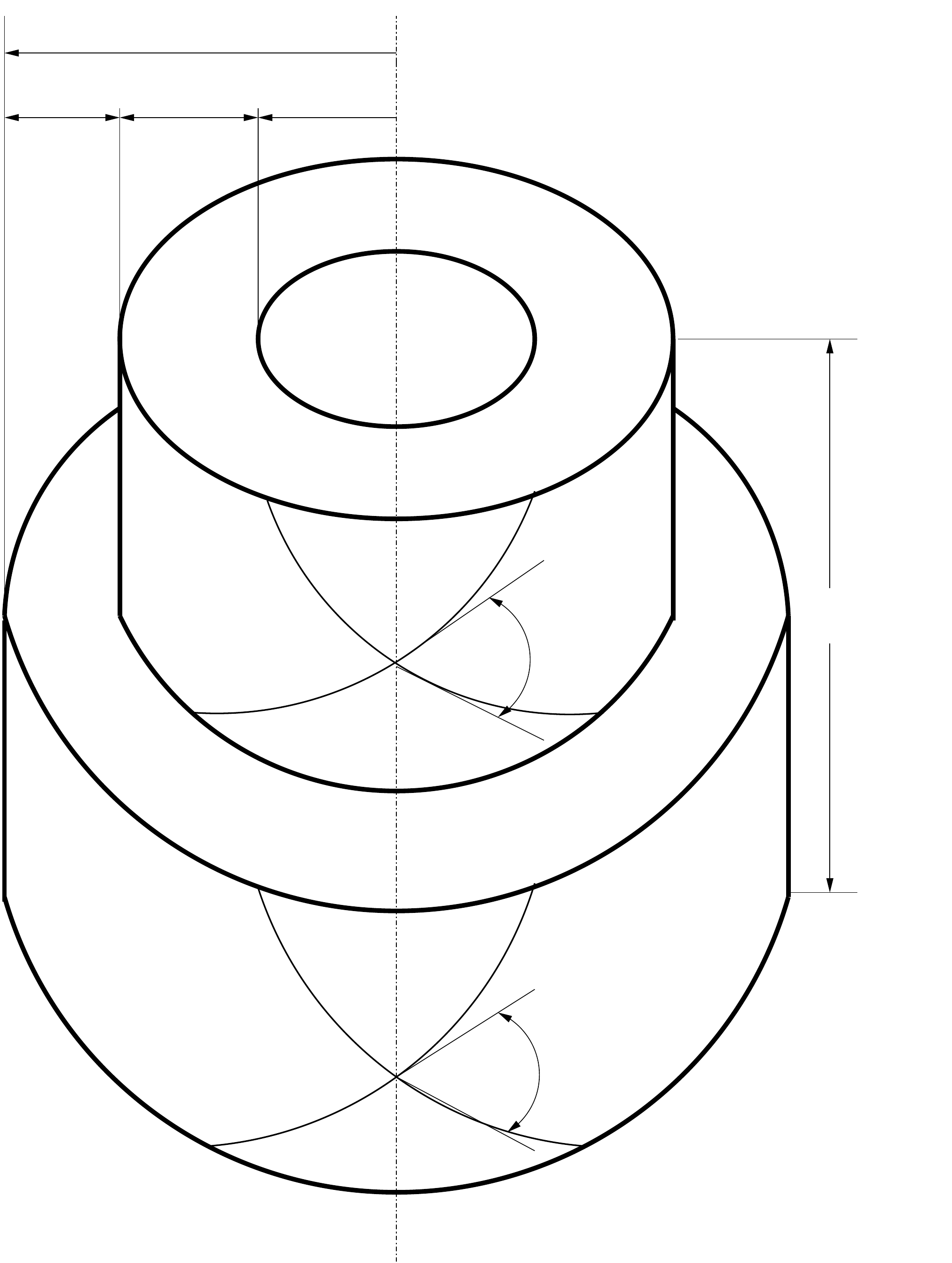}
    \caption{Geometrical data for a carotid artery from a rabbit (see \cite{Holzapfel00:00, CJC83})}\label{fig:artgeo}
  }
\end{figure}

We use the geometrical data from \cite{CJC83, Holzapfel00:00}; see an illustration 
in Fig.\ref{fig:artgeo}.
In order to compare FSI simulation using different hyperelastic models (see 
Section \ref{sec:hypermat}), we adopt the same geometrical data (except the 
angles $\alpha_M$ and $\alpha_A$) for 
the models of Neo-Hookean and Mooney-Rivlin materials. Furthermore, we 
set the value of the material parameters for three hyperelastic models as 
indicated in Tab. \ref{tab:par}, where $M$ denotes the media and $A$ the 
adventitia. 
\begin{table}[ht!]\begin{tabular}{cccccc}
\toprule
& \multicolumn{2}{c}{$c_{10}$} & \multicolumn{2}{c}{$c_{01}$} & $\rho_s$\\
 \cmidrule(r){2-3} \cmidrule(r){4-5}
& M    & A &  M    & A \\
\midrule
Neo-Hookean & $3$ kPa   & $0.3$ kPa &$-$ &$-$ & $1.2$ kg/m$^3$ \\
Mooney-Rivlin& $3$ kPa  & $0.3$ kPa &$0.3$ kPa & $0.2$ kPa  & $1.2$ kg/m$^3$ \\
Artery & $3$ kPa   & $0.3$ kPa &$-$ &$-$& $1.2$ kg/m$^3$  \\
%\bottomrule
\toprule
& \multicolumn{2}{c}{$k_1$} & \multicolumn{2}{c}{$k_2$} & $\kappa$ \\
 \cmidrule(r){2-3} \cmidrule(r){4-5}
& M    & A &  M    & A  &\\
\midrule
Neo-Hookean &$-$ &$-$ & $-$ &$-$ & $10^5$ kPa \\
Mooney-Rivlin &$-$ &$-$ & $-$ &$-$ & $10^5$ kPa \\
Artery &$2.3632$ kPa &$0.5620$ kPa & $0.8393$ &$0.7112$ & $10^5$ kPa     \\
\bottomrule
\end{tabular}\caption{The value of material parameters for three hyperelastic models.}
\label{tab:par}
\end{table}
We use Netgen \cite{JS97:00} to generate finite element mesh for 
the computational FSI domain, that provides conforming grids on the FSI interface 
and two-layered structure interface. In order to study the robustness of  
the solvers (see Section \ref{sec:lsm}) for the linearized coupled FSI system with 
respect to the discretization mesh parameter, three finite element 
meshes are generated using Netgen. In Tab. \ref{tab:mesh}, 
we summarize the total number of 
grid nodes (\#Nod), tetrahedra (\#Tet) 
and degrees of freedom (\#Dof) in the finite element simulation, that includes the mesh 
movement, fluid and structure sub-problems. 
\begin{table}[ht!]\begin{tabular}{crrr}
\toprule
&  \#Nod &  \#Tet    & \#Dof \\
\midrule
Coarse mesh & $1034$    & $4824$  & $6959$ \\
Intermediate mesh & $7249$ & $38592$  &$37909$  \\
Fine mesh & $54521$  & $308736$  &$285167$ \\
\bottomrule\end{tabular}\caption{Three finite element meshes.}
\label{tab:mesh}
\end{table}
For the fluid, we set the density 
$\rho_f=1$ mg/mm$^3$, the dynamic viscosity $\mu=0.035$ Poise. The fluid Neumann 
boundary condition on $\Gamma_{in}^t$ is given by $g_{in}=1.332 n_f$ kPa 
for $t\leq 0.125$ ms and $g_{in}=0$  kPa for $t>0.125$ ms. The remaining boundary 
conditions are specified in Section \ref{sec:pre}. The fluid and structure are at the 
rest in the initial time. The time step size $\Delta t$ is set to $0.125$ ms. We run 
the simulation until $12$ ms. 
\subsection{Convergence of Newton's method}\label{sec:pnewt}
To verify the linearization for the coupled nonlinear FSI system 
(see Section \ref{sec:nsm}), we show the relative error (err) and 
iteration number (\#it) of Newton's method for the FSI simulation using three 
different hyperelastic models: Neo-Hookean (FSI\_NH), Mooney-Rivlin (FSI\_MR) and 
artery (FSI\_AR), and three different meshes: Coarse mesh (C), 
intermediate mesh (I) and fine mesh (F); see Tab. \ref{tab:pnewt} for details. Note 
that since we observe the same performance of Newton's method for solving 
the nonlinear system at all time steps, only the performance at the first time 
step is recorded in Tab. \ref{tab:pnewt} for simplicity of presentation.  
\begin{table}[ht!]\begin{tabular}{cccc}
\toprule
\multicolumn{1}{l} {FSI\_NH:} \\ 
\cmidrule(r){1-1} 
\#it &\multicolumn{3}{c}{err}  \\
 \cmidrule(r){2-4} 
& C    & I & F    \\
\midrule
1 &$6.2e+01$ &$6.8e+01$ &$6.9e+01$ \\
2 &$3.8e-02$  &$5.7e-02$  &$6.6e-02$\\
3 &$7.4e-06$  &$7.7e-06$  &$4.0e-06$\\
4 &$4.3e-09$  &$1.1e-09$  &$1.3e-08$\\
\midrule
\multicolumn{1}{l} {FSI\_MR:}\\
 \cmidrule(r){1-1} 
 \#it  & \multicolumn{3}{c}{err} \\
 \cmidrule(r){2-4} 
& C    & I & F  \\
\midrule
1& $6.2e+01$  &$6.8e+01$ &$6.9e+01$\\
2& $3.8e-02$   &$5.7e-02$  &$6.4e-02$\\
3& $7.3e-06$   &$7.3e-06$  &$3.0e-06$\\
4& $4.3e-09$   &$9.7e-10$  &$2.5e-09$\\
\midrule
\multicolumn{1}{l} {FSI\_AR:}\\
 \cmidrule(r){1-1} 
 \#it & \multicolumn{3}{c}{err}  \\
 \cmidrule(r){2-4} 
& C    & I & F \\
\midrule
1& $6.2e+01$ &$6.8e+01$ &$6.9e+01$\\
2& $4.0e-02$  &$6.0e-02$  &$7.0e-02$ \\
3& $7.5e-06$  &$8.4e-06$  &$1.2e-05$ \\
4& $4.3e-09$  &$1.2e-09$  &$4.0e-09$ \\
\bottomrule
\end{tabular}\caption{The convergence history of Newton's method for the FSI nonlinear system using three hyperelastic models and meshes.}
\label{tab:pnewt}
\end{table}
From the convergence history displayed in Tab. \ref{tab:pnewt}, we observe (near)quadratic convergence rate of Newton's method, that 
conforms the derivation for the linearization of the coupled nonlinear FSI system, stemming from the domain movements, 
convection terms, material laws, transmission conditions and stabilization parameters. We observe 
nearly the same convergence rate for the nonlinear FSI system using three different hyperelastic models 
on the coarse, intermediate and fine mesh. At each iteration of Newton's method, we use the 
preconditioned Krylov subspace, algebraic multigrid and multilevel methods to solve the linearized FSI system; 
see numerical results in Section \ref{sec:plins} and \ref{sec:amlins}. 
\subsection{Iteration numbers of preconditioned Krylov subspace methods}\label{sec:plins}
To compare performance of preconditioned Krylov subspace methods 
for the linearized coupled FSI system, we use the GMRES method combined 
with the preconditioners from Section \ref{sec:pky}. The stopping criterion for the 
GMRES method is set by the relative error $10^{-9}$. We compare the total number 
of GMRES iterations (\#it) 
to reach this criterion for the FSI simulation using 
the Neo-Hookean (FSI\_NH), Mooney-Rivlin (FSI\_MR) and 
artery (FSI\_AR) model on coarse mesh (C), 
intermediate mesh (I) and fine mesh (F). The detailed numerical 
results are shown in Tab. \ref{tab:pkrylov}. Note that since the performance is 
similar for all Newton iterations, we demonstrate the iteration numbers 
at the first Newton iteration. The inverse of each sub-problem in 
the preconditioners is realized by calling the corresponding AMG cycle, 
that has been developed in \cite{ULHY13}.   
\begin{table}[ht!]\begin{tabular}{cccccccccc}
\toprule
Precontitioner &\multicolumn{9}{c}{\#it}  \\
 \cmidrule(r){2-10} 
 &\multicolumn{3}{c} {FSI\_NH} &\multicolumn{3}{c} {FSI\_MR}
&\multicolumn{3}{c} {FSI\_AR} \\ 
\cmidrule(r){2-4} \cmidrule(r){5-7}  \cmidrule(r){8-10} 
& C & I & F & C & I & F & C & I & F \\
\midrule
$\tilde{P}_D$         &$51$ & $111$  &$217$  &$53$  &$111$  &$227$       &$46$ &$98$ & $189$\\
$\tilde{P}_L$         &$28$ & $58$    & $109$ &$29$  &$60$    &$114$   &$25$ & $50$ &$95$  \\
$\tilde{P}_U$         &$28$ & $59$   & $114$ &$28$  &$61$     &$119$  &$25$ &$51$ &$98$  \\
$\tilde{P}_{SSOR}$ &$27$ & $54$   & $104$ &$28$  &$57$    & $108$  &$24$ &$48$ &$91$  \\
$\tilde{P}_{ILU}$   &$27$ & $54$    & $104$ &$28$  &$57 $   &$108$   &$24$ &$48$ &$91$  \\
\bottomrule
\end{tabular}\caption{The performance of preconditioned GMRES method 
    for the linearized FSI system using three hyperelastic models and meshes.}
\label{tab:pkrylov}
\end{table}
As we observe from the iteration numbers of the linear solvers using different 
preconditioners in Tab. \ref{tab:pkrylov}, the solver with the preconditioner $\tilde{P}_D$ 
requires more iteration numbers than the other four preconditioners. The solvers 
with the preconditioners $\tilde{P}_L$, $\tilde{P}_U$, $\tilde{P}_{SSOR}$ and $\tilde{P}_{ILU}$ require 
almost the same number of iteration numbers. As expected, when the mesh is refined, 
the iteration number of the preconditioned GMRES method increases. We will see 
in Section \ref{sec:amlins} that, the mesh dependence is eliminated by using 
the multigrid and multilevel method. 
\subsection{Iteration numbers of algebraic multigrid and multilevel methods}
\label{sec:amlins}
In this section, we compare the performance of the AMG and AMLI 
method for the linearized coupled FSI system. More precisely, 
we show the number of iteration numbers (\#it) of the AMG, AMLI, 
AMG preconditioned GMRES (AMG\_GMRES), AMG preconditioned 
FGMRES (AMG\_FGMRES),  
AMLI preconditioned GMRES (AMLI\_GMRES) and 
AMLI preconditioned FGMRES (AMLI\_FGMRES) method, respectively, 
up to the relative error $10^{-9}$. 
We run the FSI simulation using the Neo-Hookean (FSI\_NH), 
Mooney-Rivlin (FSI\_MR), artery (FSI\_AR) model, on the coarse (C), intermediate (I) 
and fine (F) mesh, respectively. See Tab. \ref{tab:pkamg} for details. 
We use $8-10$ smoothing steps in the AMG and AMLI cycle, each of which only 
requires $1$ AMG cycle for the corresponding mesh movement, 
fluid and structure sub-problem (see Section \ref{sec:fsism}). As preconditioners, 
we only apply $1$ AMG or AMLI cycle in the preconditioned 
GMRES or FGMRES iteration.
\begin{table}[ht!]\begin{tabular}{cccccccccc}
\toprule
Method &\multicolumn{9}{c}{\#it}  \\
 \cmidrule(r){2-10} 
 &\multicolumn{3}{c} {FSI\_NH} &\multicolumn{3}{c} {FSI\_MR}
&\multicolumn{3}{c} {FSI\_AR} \\ 
\cmidrule(r){2-4} \cmidrule(r){5-7}  \cmidrule(r){8-10} 
& C    & I & F    & C    & I & F & C    & I & F \\
\midrule
AMG                     &$7$ &$8$ & $12$   &$8$ &$8$ & $11$    &$7$ &$8$  & $10$ \\
AMG\_GMRES     &$7$ &$7$ & $8$   &$7$ &$7$ & $8$   &$7$ &$7$ & $8$ \\
AMG\_FGMRES   &$6$ &$7$ & $8$   &$6$ &$7$ & $9$   &$6$ &$7$ & $8$ \\
AMLI                    &$7$ &$8$ & $12$ &$8$ &$8$ & $11$  &$7$ &$8$ & $10$ \\
AMLI\_GMRES    &$7$ &$7$ & $8$   &$7$ &$7$ & $8$    &$7$ &$7$ & $8$ \\
AMLI\_FGMRES  &$6$ &$7$ & $8$   &$6$ &$7$ & $9$    &$6$ &$7$ & $8$\\
\bottomrule
\end{tabular}\caption{The performance of the AMG, AMLI, and AMG and 
      AMLI preconditioned Krylov subspace method for the linearized FSI
      system using three hyperelastic models and meshes.}
\label{tab:pkamg}
\end{table}
As we observe from Tab. \ref{tab:pkamg}, the AMG and AMLI method 
requires the same iteration numbers for each case. The AMG and AMLI 
preconditioned GMRES and FGMRES methods show improved performance 
with fewer iteration numbers than the AMG and AMLI methods. When the 
mesh is refined, we observe the iteration numbers using these methods stay 
in a very similar range. This demonstrates the robustness of the multigrid 
and multilevel method for the coupled FSI system with respect to the mesh 
refinement. 

\subsection{Visualization of  the numerical solutions}\label{sec:vis}
\begin{figure}[htbp]
  \centering{
    \begin{subfigure}[b]{0.32\textwidth}
      \includegraphics[scale=0.32]{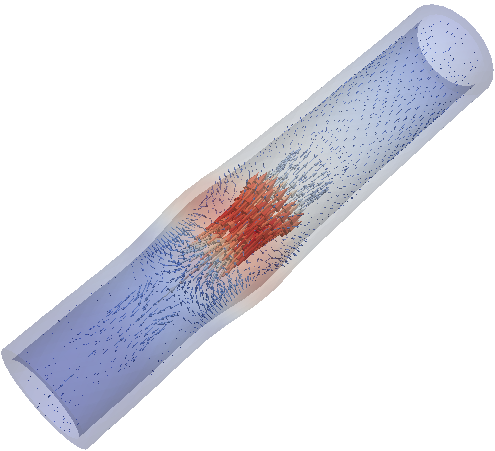}
      \caption{Neo-Hookean}
    \end{subfigure}
    \begin{subfigure}[b]{0.32\textwidth}
      \includegraphics[scale=0.32]{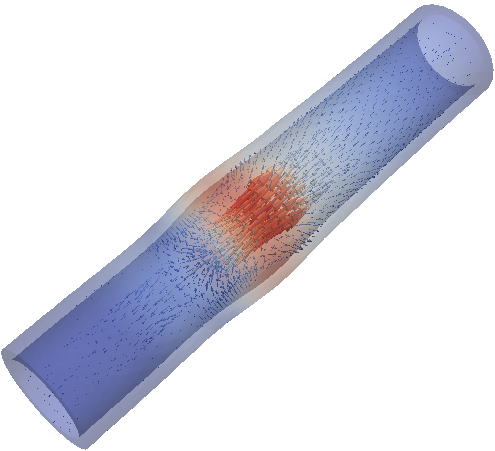}
      \caption{Mooney-Rivlin}
    \end{subfigure}
    \begin{subfigure}[b]{0.32\textwidth}
      \includegraphics[scale=0.32]{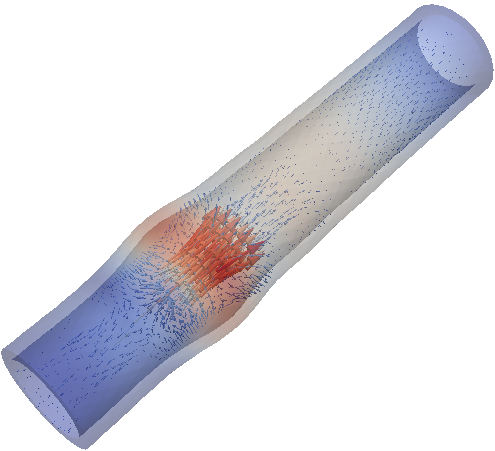}
      \caption{Artery}
    \end{subfigure}
  }
  \caption{Visualization of structure deformation and fluid velocity fields at time 
    $t=8$ ms using 
    three different structure models of Neo-Hookean material (left), 
    Mooney-Rivlin material (middle) 
    and anisotropic two-layer thick walled artery.}
  \label{fig:numsol}
\end{figure}

In order to demonstrate the numerical simulation results, 
we visualize the structure deformations and fluid velocity fields in Fig. \ref{fig:numsol}, 
where the FSI solutions at time level $t=8$ ms using the structure models of 
the Neo-Hookean material, the Mooney-Rivlin material and the 
anisotropic two-layer thick walled artery, are respectively shown. 
\subsection{Comparison with the partitioned approach}\label{sec:comp}
In this section, we compare the numerical simulation results 
obtained by the monolithic approach with the results by the 
partitioned approach as in \cite{ULHY13}. 

\begin{figure}[htbp]
  \centering{
    \begin{subfigure}[b]{0.4\textwidth}
      \includegraphics[scale=0.125]{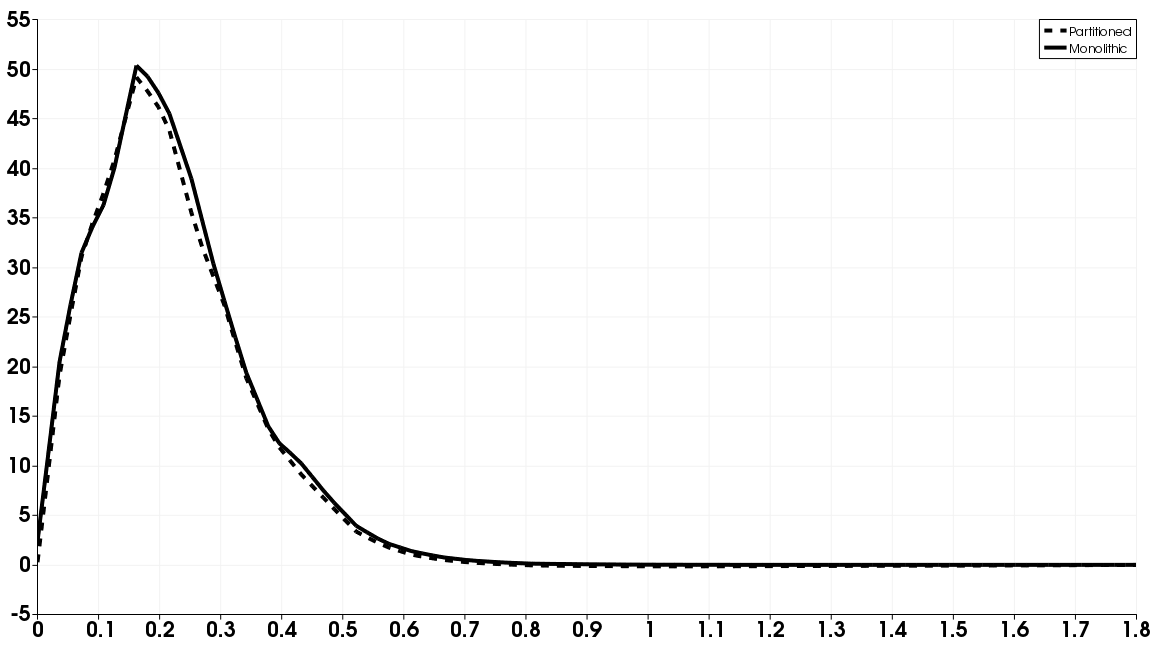}
      \caption{$t=1.50$ ms}
    \end{subfigure}
    \hfill
    \begin{subfigure}[b]{0.4\textwidth}
      \includegraphics[scale=0.125]{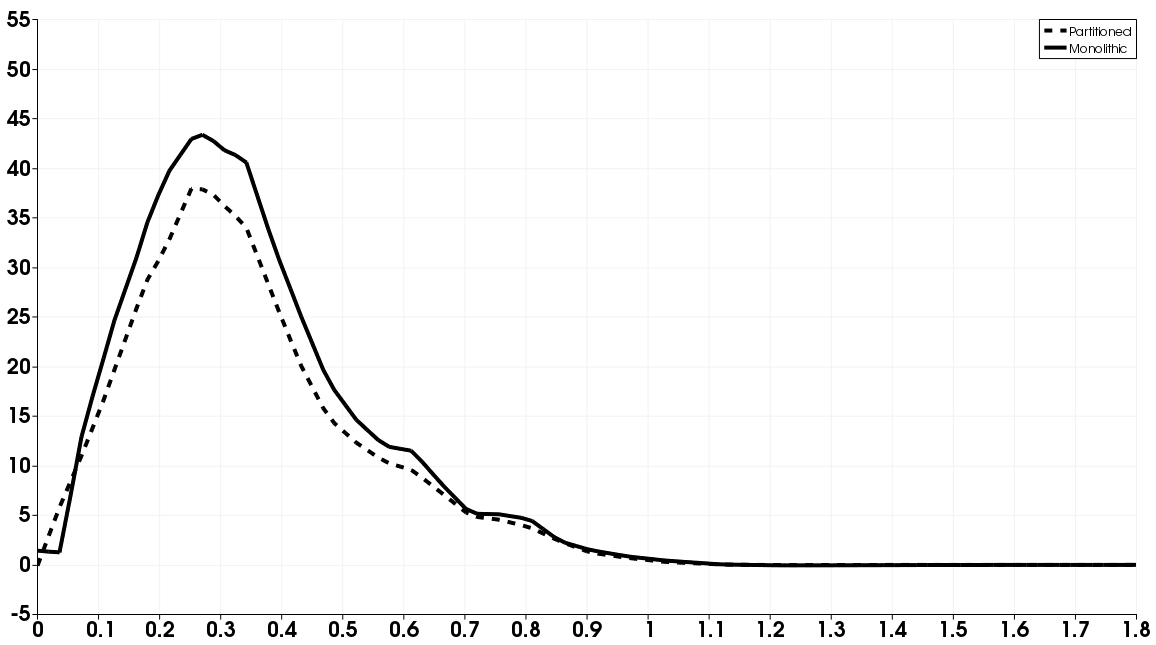}
      \caption{$t=3.00$ ms}
    \end{subfigure}
    \begin{subfigure}[b]{0.4\textwidth}
      \includegraphics[scale=0.125]{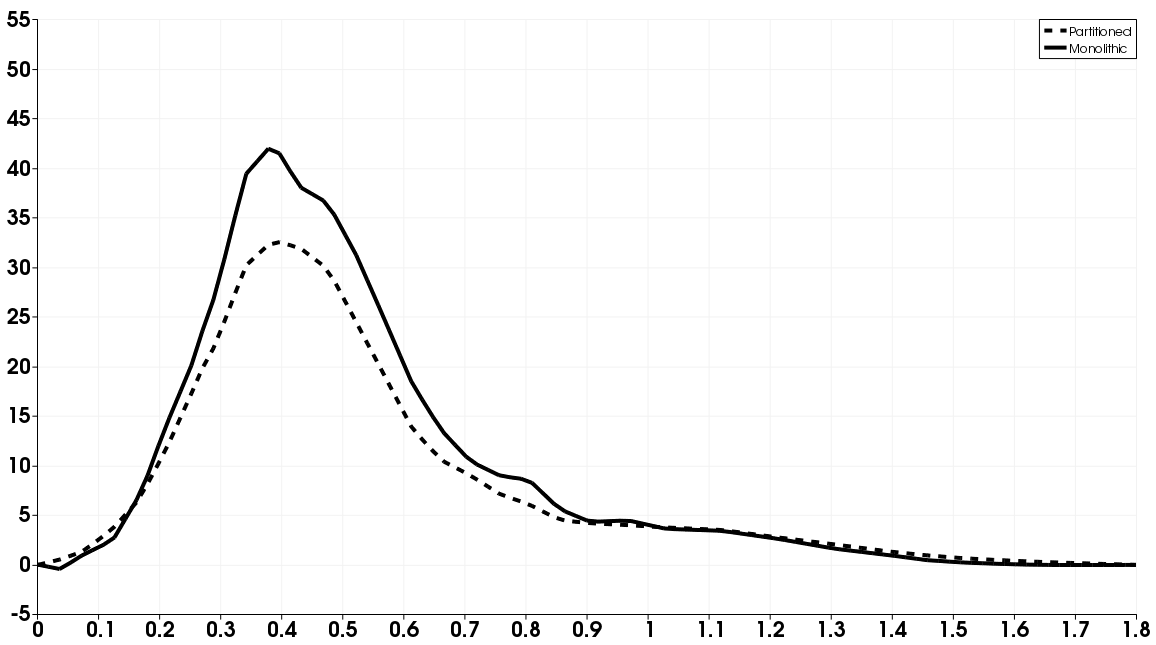}
      \caption{$t=4.50$ ms}
    \end{subfigure}
    \hfill
    \begin{subfigure}[b]{0.4\textwidth}
      \includegraphics[scale=0.125]{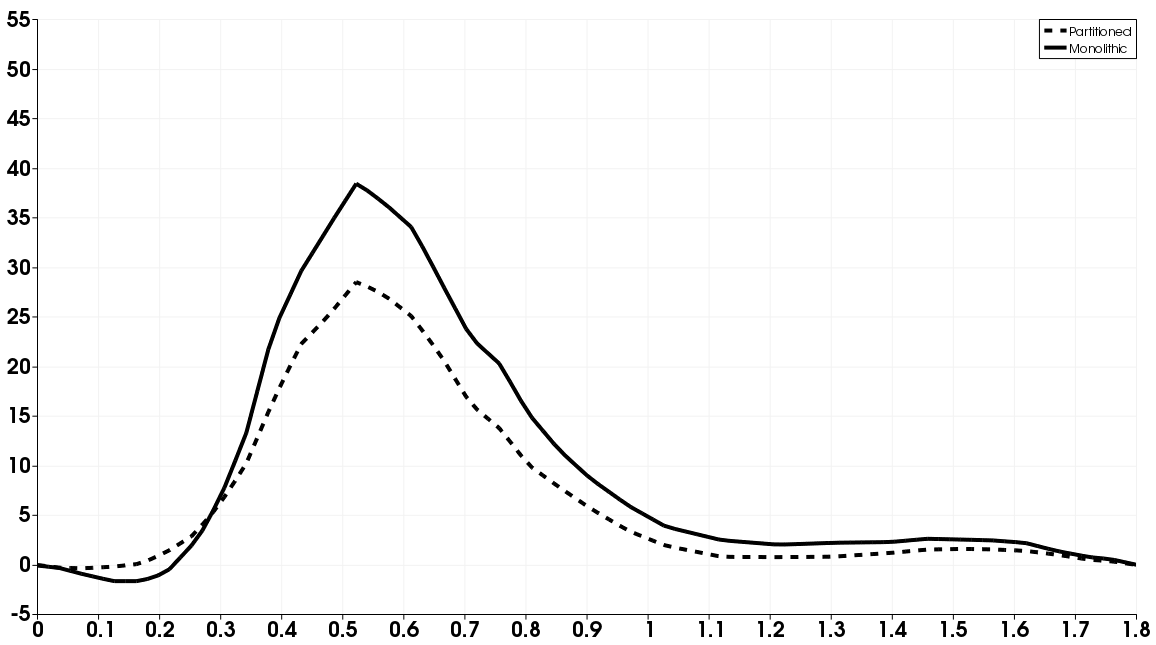}
      \caption{$t=6.00$ ms}
    \end{subfigure}
    \begin{subfigure}[b]{0.4\textwidth}
      \includegraphics[scale=0.125]{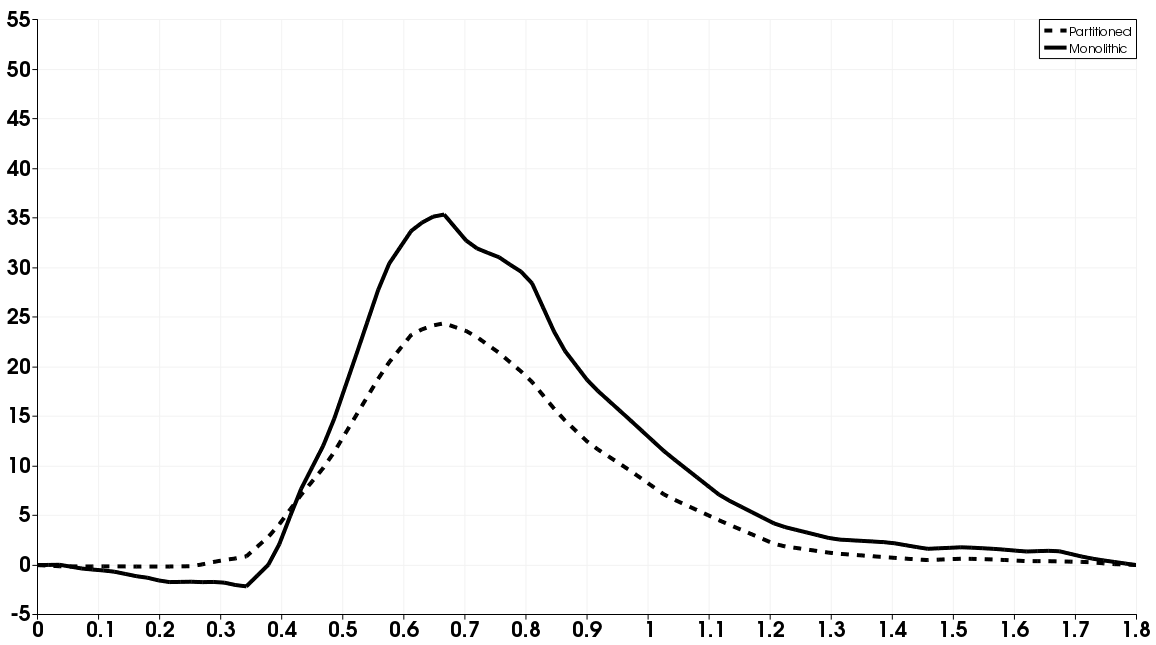}
      \caption{$t=7.50$ ms}
    \end{subfigure}
    \hfill
    \begin{subfigure}[b]{0.4\textwidth}
      \includegraphics[scale=0.125]{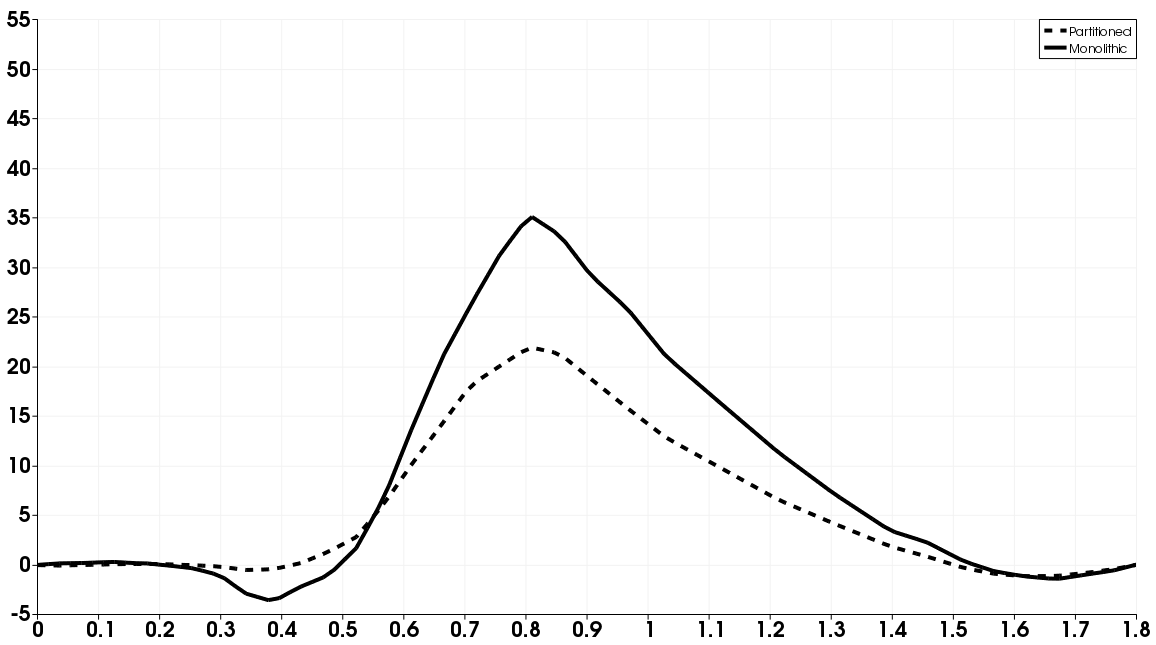}
      \caption{$t=9.00$ ms}
    \end{subfigure}
  }
  \caption{Comparison of fluid pressure waves  along the center line with the 
    starting point $(0, 0, 0)$ cm and ending point $(0, 0, 1.8)$ cm 
    from the FSI simulation using the model of Neo-Hookean material: 
    Monolithic solution (in solid lines) and partitioned solution (in dashed lines).}
  \label{fig:cmp_press_nh}
\end{figure}

\begin{figure}[htbp]
  \centering{
    \begin{subfigure}[b]{0.4\textwidth}
      \includegraphics[scale=0.125]{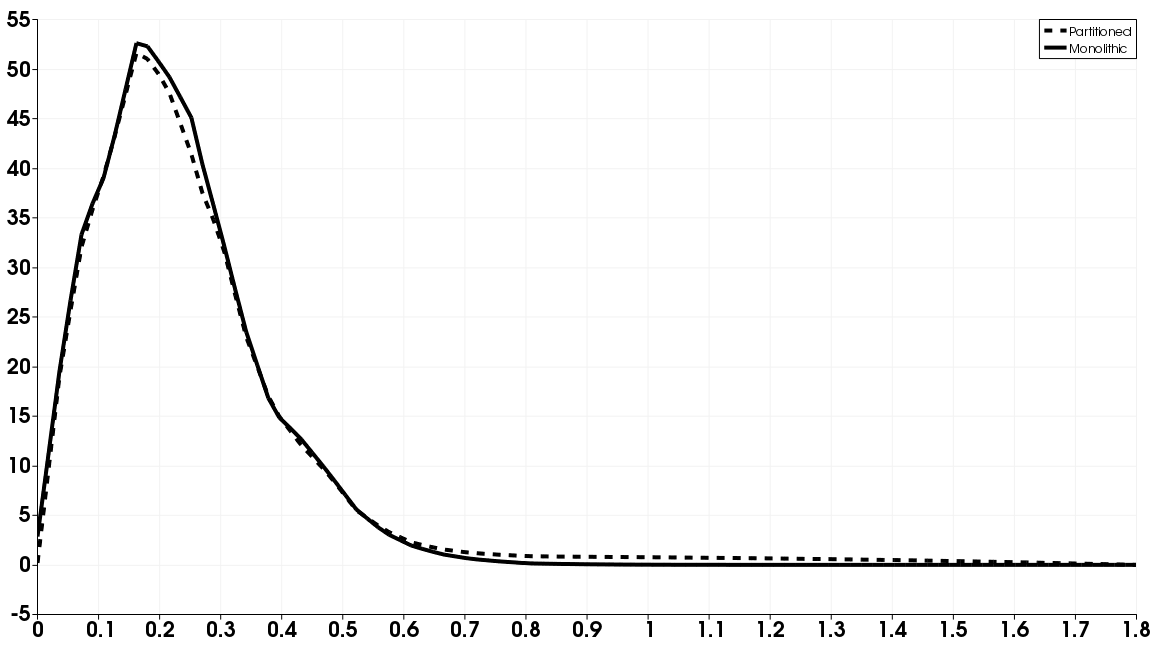}
      \caption{$t=1.50$ ms}
    \end{subfigure}
    \hfill
    \begin{subfigure}[b]{0.4\textwidth}
      \includegraphics[scale=0.125]{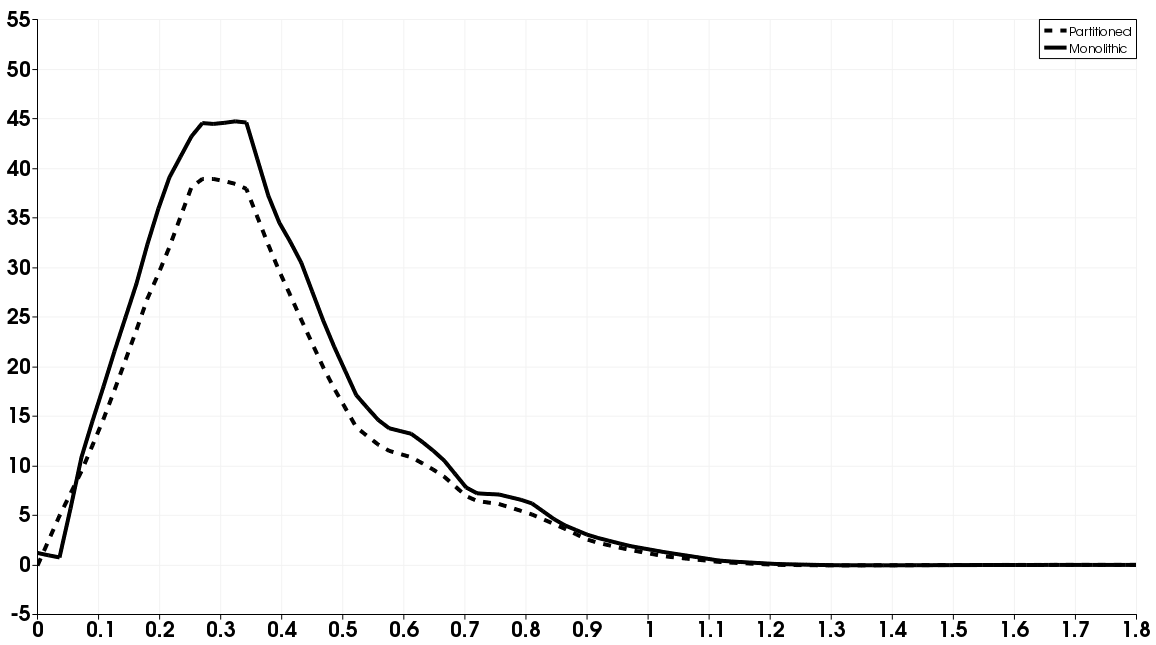}
      \caption{$t=3.00$ ms}
    \end{subfigure}
    \begin{subfigure}[b]{0.4\textwidth}
      \includegraphics[scale=0.125]{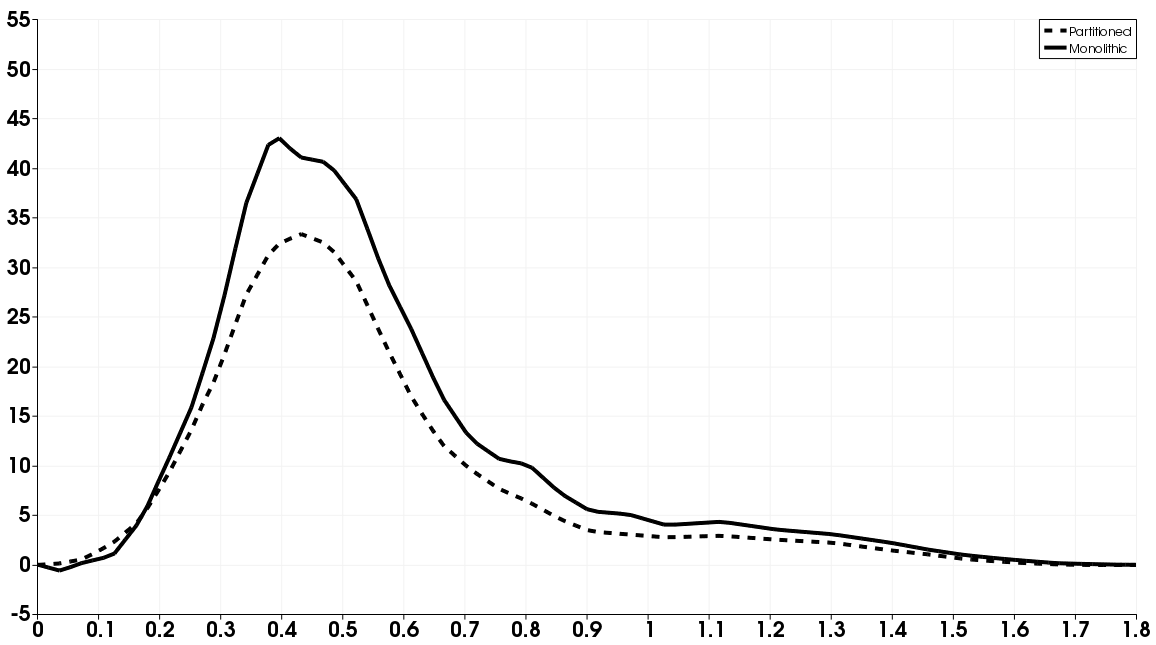}
      \caption{$t=4.50$ ms}
    \end{subfigure}
    \hfill
    \begin{subfigure}[b]{0.4\textwidth}
      \includegraphics[scale=0.125]{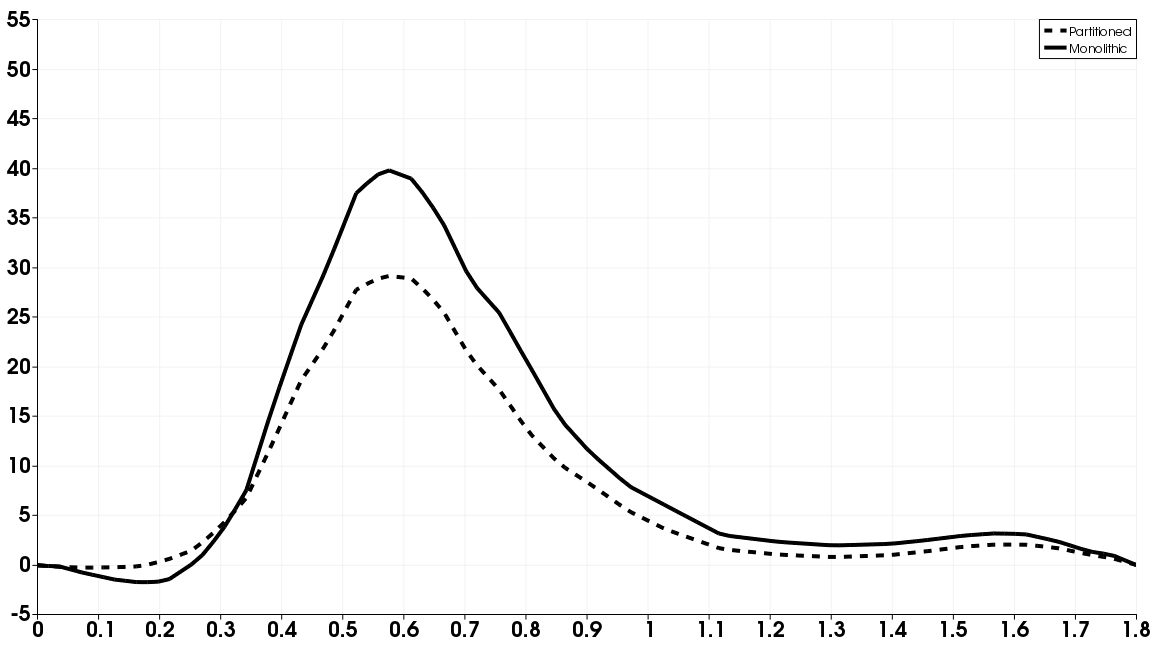}
      \caption{$t=6.00$ ms}
    \end{subfigure}
    \begin{subfigure}[b]{0.4\textwidth}
      \includegraphics[scale=0.125]{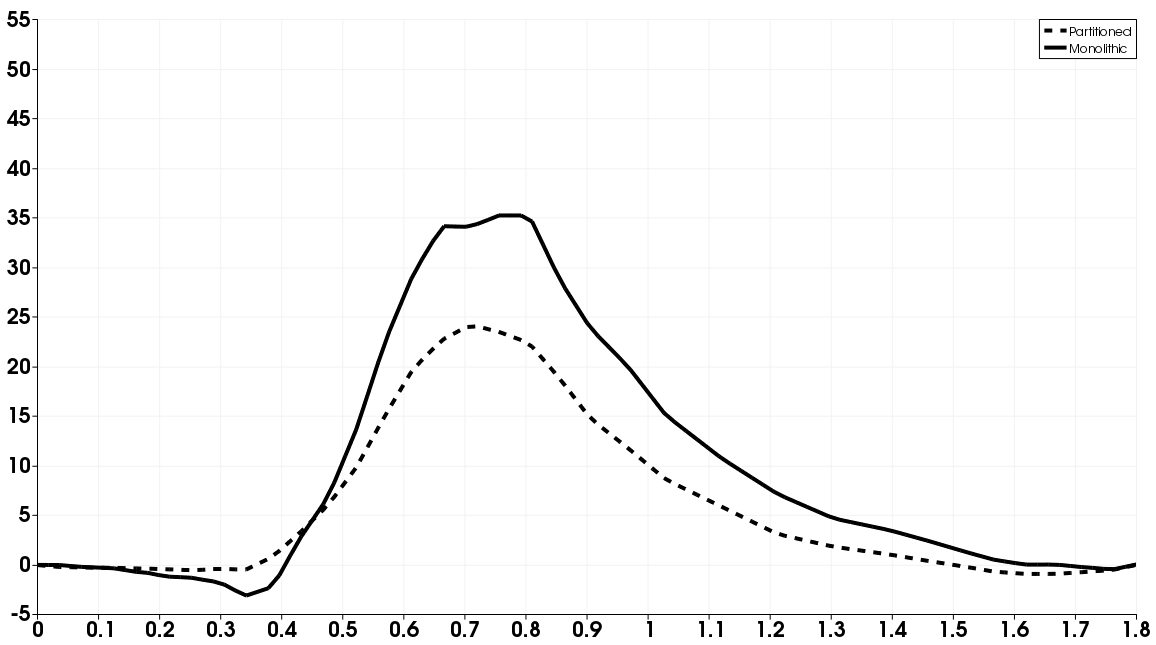}
      \caption{$t=7.50$ ms}
    \end{subfigure}
    \hfill
    \begin{subfigure}[b]{0.4\textwidth}
      \includegraphics[scale=0.125]{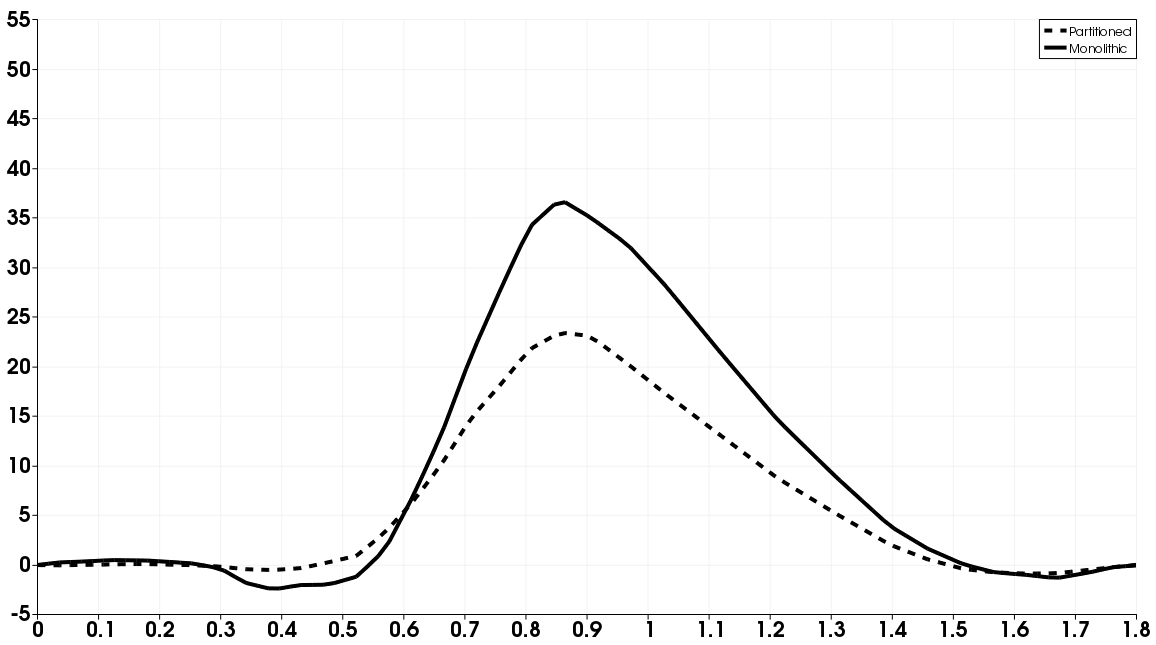}
      \caption{$t=9.00$ ms}
    \end{subfigure}
  }
  \caption{Comparison of fluid pressure waves  along the center line with the 
    starting point $(0, 0, 0)$ cm and ending point $(0, 0, 1.8)$ cm 
    from the FSI simulation using the modified model of 
    Mooney-Rivlin material: 
    Monolithic solution (in solid lines) and partitioned solution (in dashed lines).}
  \label{fig:cmp_press_mr}
\end{figure}

\begin{figure}[htbp]
  \centering{
    \begin{subfigure}[b]{0.4\textwidth}
      \includegraphics[scale=0.125]{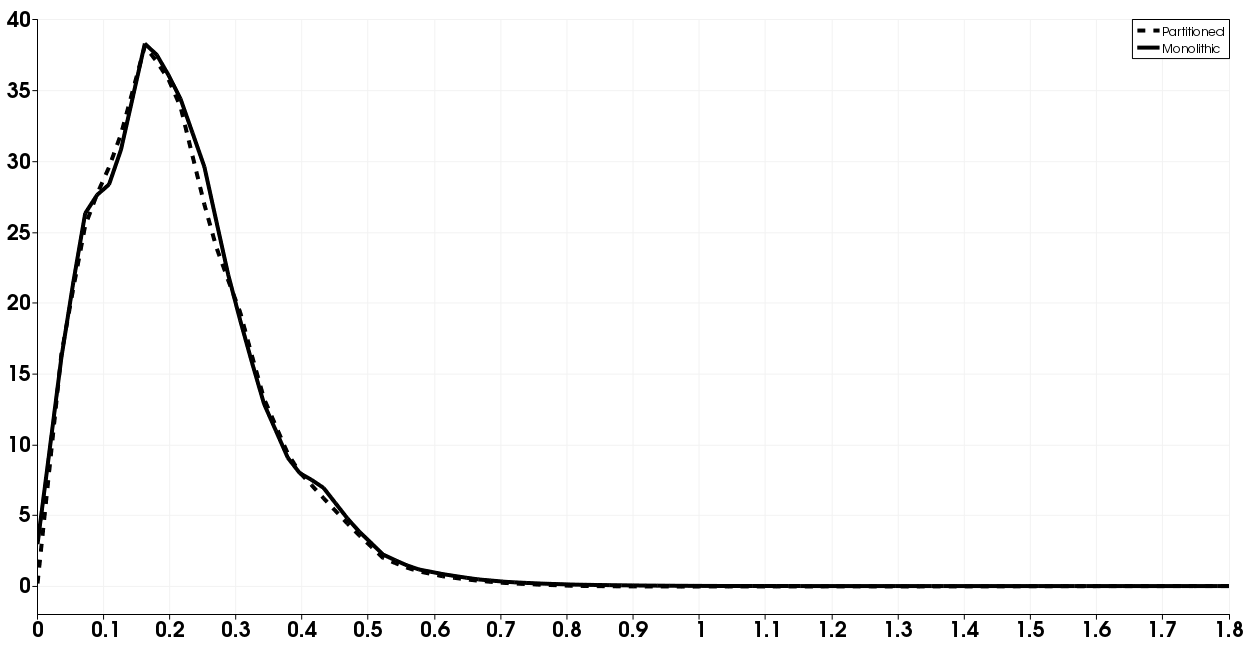}
      \caption{$t=1.75$ ms}
    \end{subfigure}
    \hfill
    \begin{subfigure}[b]{0.4\textwidth}
      \includegraphics[scale=0.125]{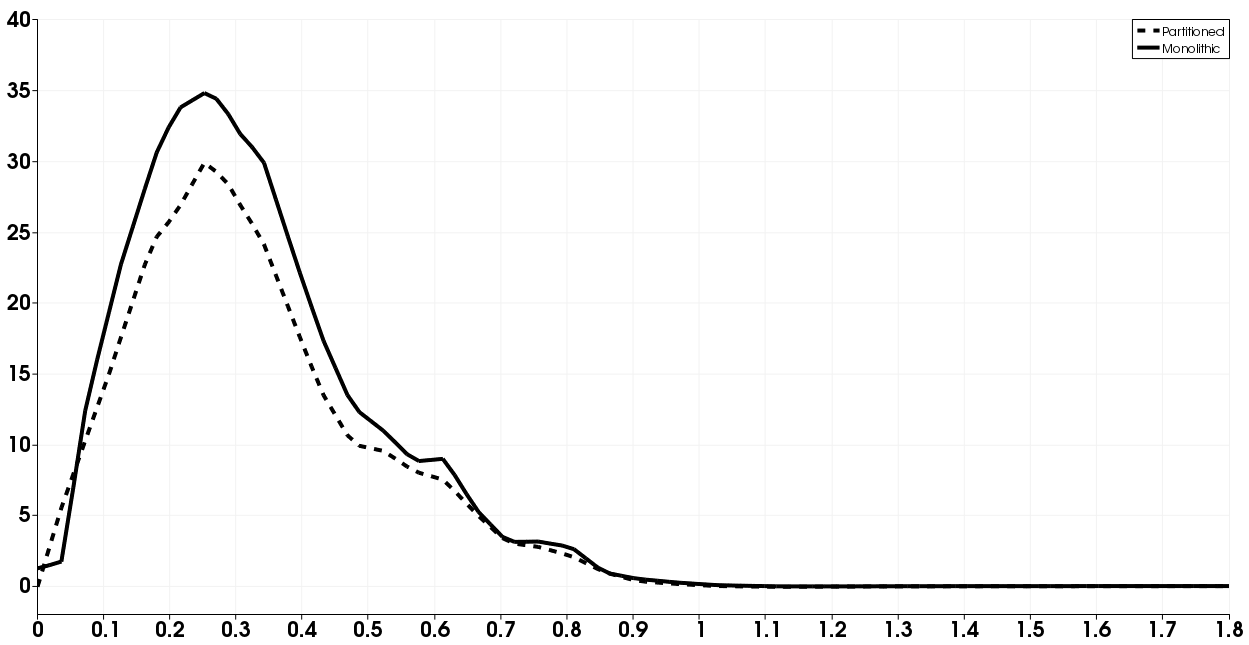}
      \caption{$t=3.50$ ms}
    \end{subfigure}
    \begin{subfigure}[b]{0.4\textwidth}
      \includegraphics[scale=0.125]{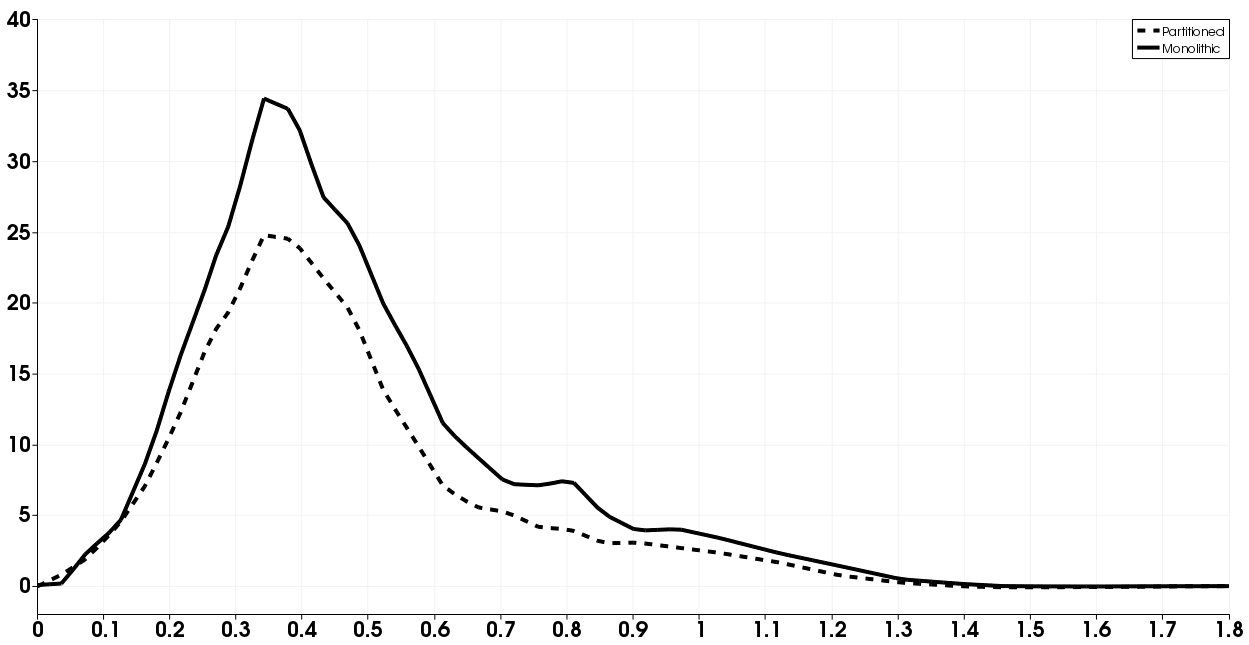}
      \caption{$t=5.25$ ms}
    \end{subfigure}
    \hfill
    \begin{subfigure}[b]{0.4\textwidth}
      \includegraphics[scale=0.125]{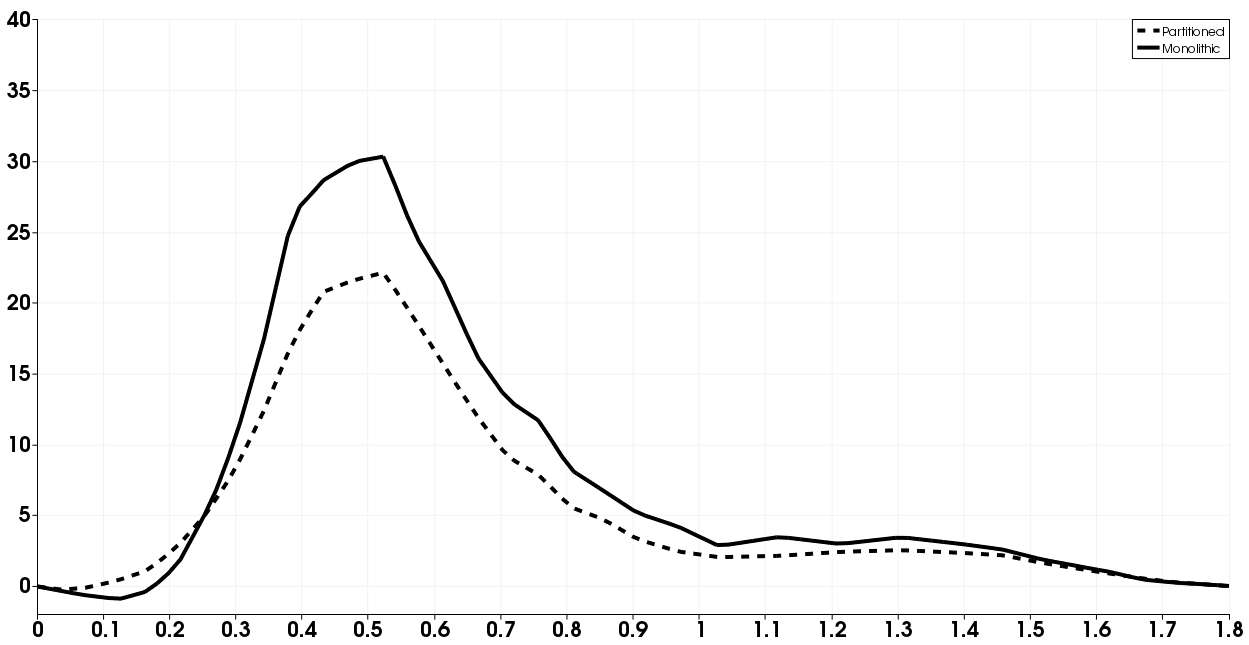}
      \caption{$t=7.00$ ms}
    \end{subfigure}
    \begin{subfigure}[b]{0.4\textwidth}
      \includegraphics[scale=0.125]{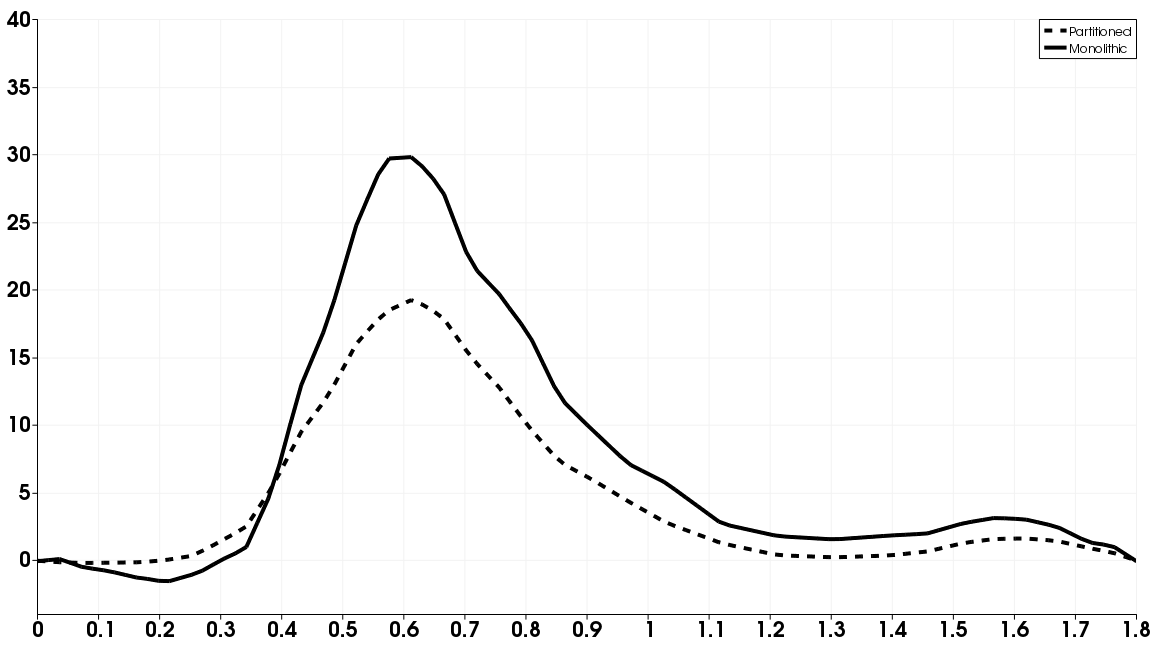}
      \caption{$t=8.75$ ms}
    \end{subfigure}
    \hfill
    \begin{subfigure}[b]{0.4\textwidth}
      \includegraphics[scale=0.125]{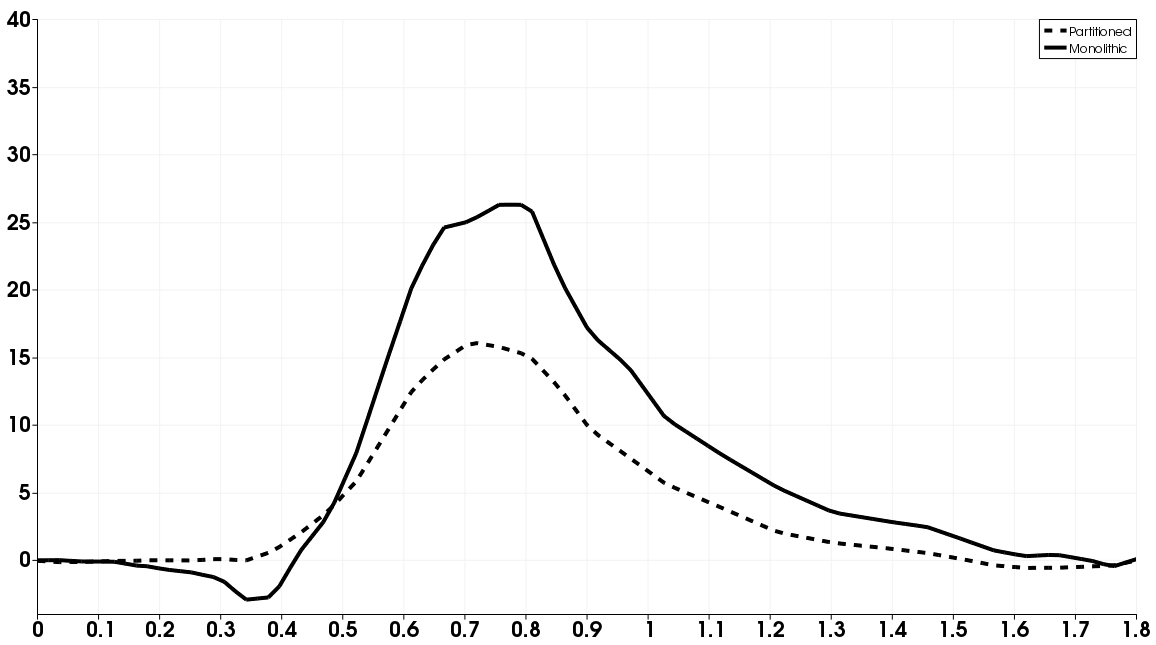}
      \caption{$t=10.50$ ms}
    \end{subfigure}
  }
  \caption{Comparison of fluid pressure waves  along the center line with the 
    starting point $(0, 0, 0)$ cm and ending point $(0, 0, 1.8)$ cm 
    from the FSI simulation using the anisotropic two-layer thick walled artery: 
    Monolithic solution (in solid lines) and partitioned solution (in dashed lines).}
  \label{fig:cmp_press_ga}
\end{figure}

We first compare the fluid pressure waves obtained from the 
FSI simulation using different structure models. 
In Fig. \ref{fig:cmp_press_nh}, Fig. \ref{fig:cmp_press_mr} and 
Fig. \ref{fig:cmp_press_ga}, we plot fluid pressure waves along the center line 
with the starting point $(0, 0, 0)$ cm and ending point $(0, 0, 1.8)$ cm, for the 
model of Neo-Hookean material, Mooney-Rivlin material and 
anisotropic two-layer thick walled artery, respectively. In each subplot of these 
three figures, the horizontal line represents the center line (in cm), and the 
vertical line represents the pressure (in Pa). 

We compare the pressure 
waves at different time levels using the monolithic and partitioned approach. 
According to our experiments, we observe at the first time steps, the solution 
obtained by using the monolithic and partitioned approach conforms to each 
other very well. With time stepping, the solution obtained by the partitioned 
approach has smaller magnitude than the solution by the monolithic approach. 
This is due to the fact that, 
at each time level of the partitioned approach,
%in the partitioned approach, at each time level 
we apply the fixed-point method 
to the reduced interface equation in an iterative manner, which introduce 
some additional errors in the solution procedure. These additional 
errors are accumulated with time stepping. However, for the 
monolithic approach, we solve the coupled system in an all-at-once manner, 
such additional errors are eliminated. 

\begin{figure}[htbp]
   \centering{
    \begin{subfigure}[b]{0.4\textwidth}
      \includegraphics[scale=0.125]{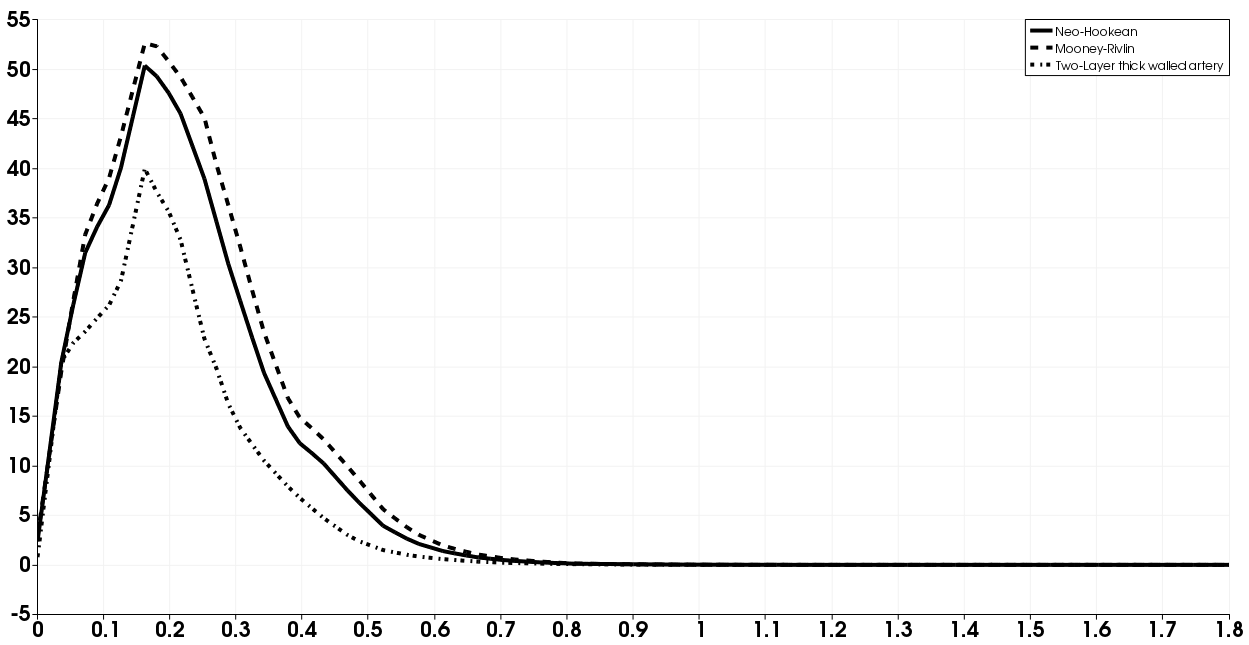}
      \caption{$t=1.50$ ms}
    \end{subfigure}
    \hfill
    \begin{subfigure}[b]{0.4\textwidth}
      \includegraphics[scale=0.125]{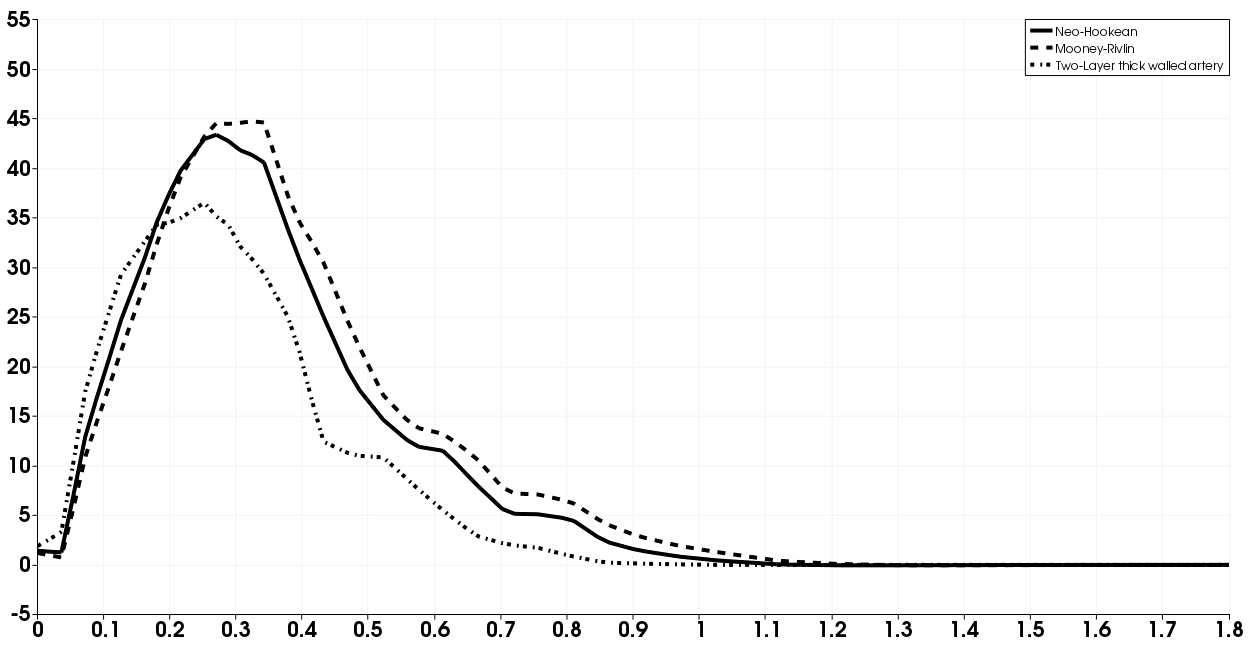}
      \caption{$t=3.00$ ms}
    \end{subfigure}
    \begin{subfigure}[b]{0.4\textwidth}
      \includegraphics[scale=0.125]{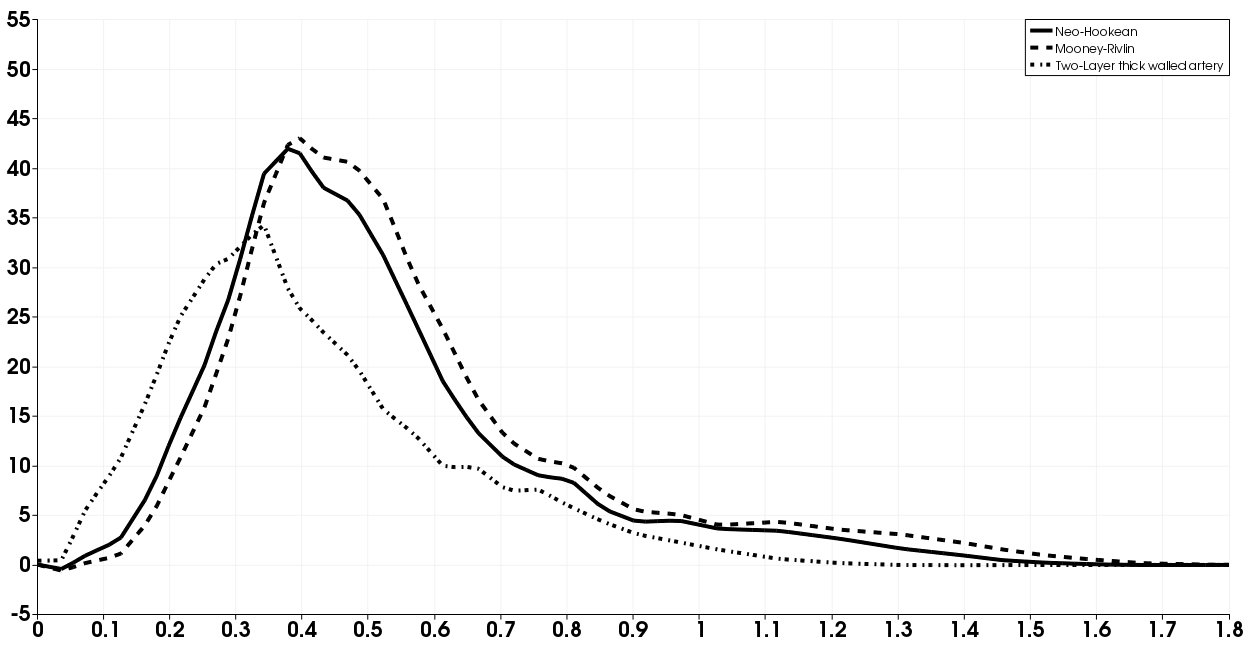}
      \caption{$t=4.50$ ms}
    \end{subfigure}
    \hfill
    \begin{subfigure}[b]{0.4\textwidth}
      \includegraphics[scale=0.125]{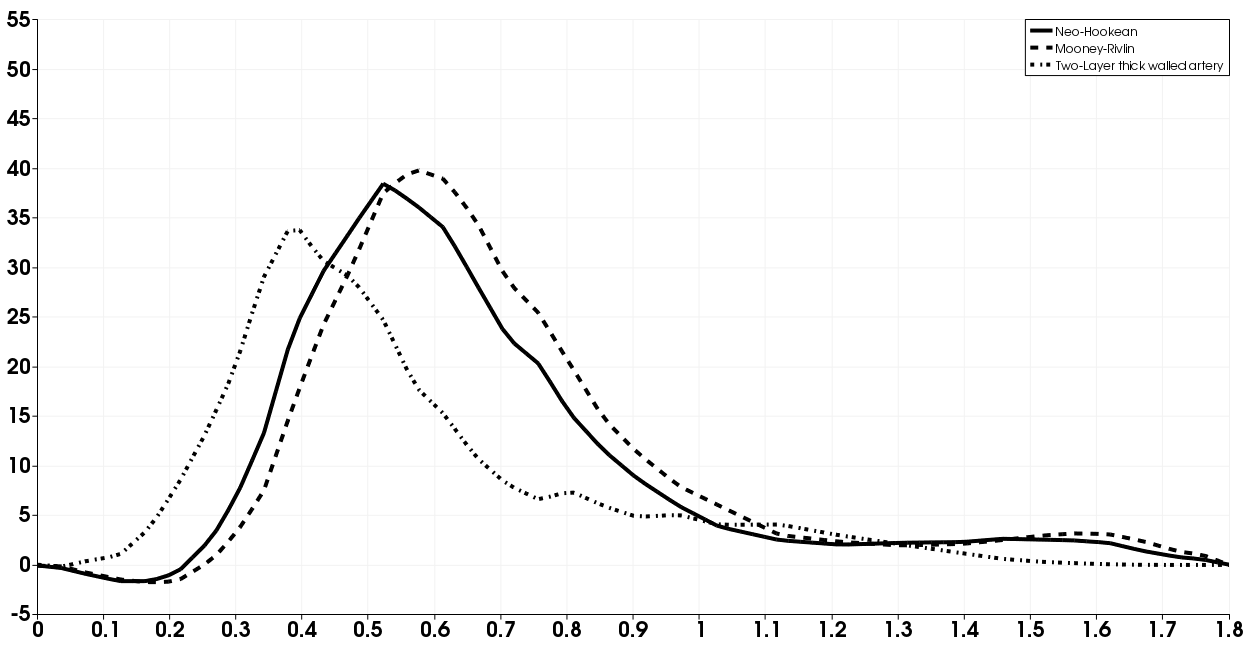}
      \caption{$t=6.00$ ms}
    \end{subfigure}
    \begin{subfigure}[b]{0.4\textwidth}
      \includegraphics[scale=0.125]{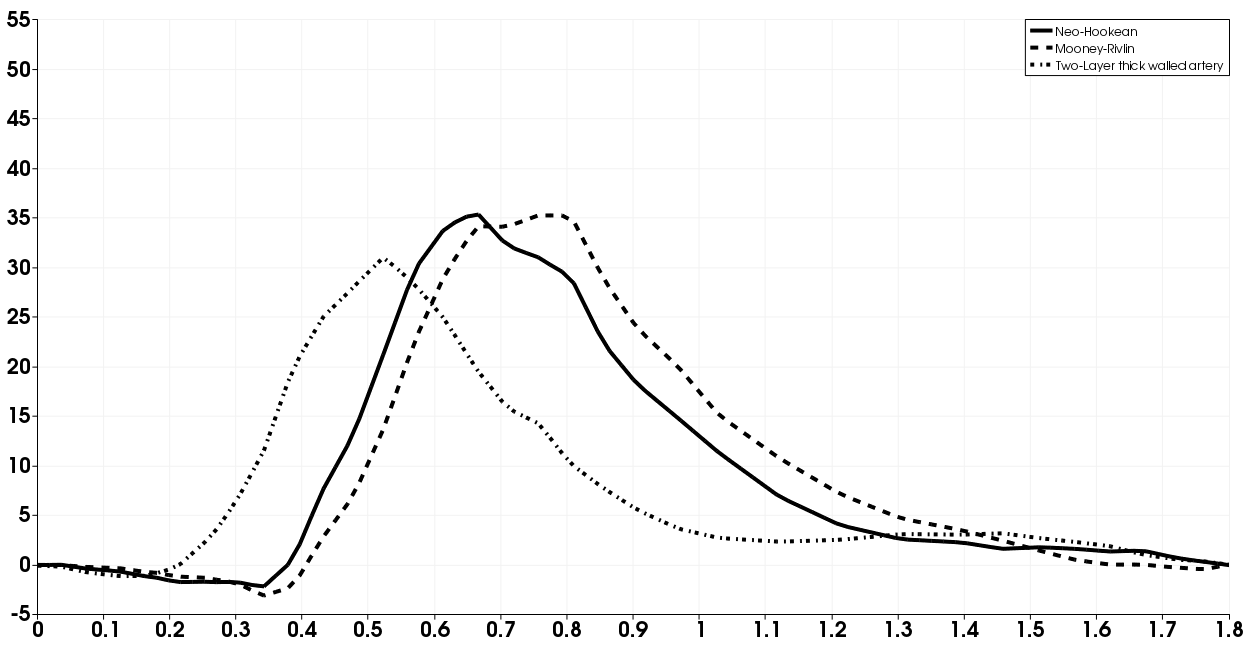}
      \caption{$t=7.50$ ms}
    \end{subfigure}
    \hfill
    \begin{subfigure}[b]{0.4\textwidth}
      \includegraphics[scale=0.125]{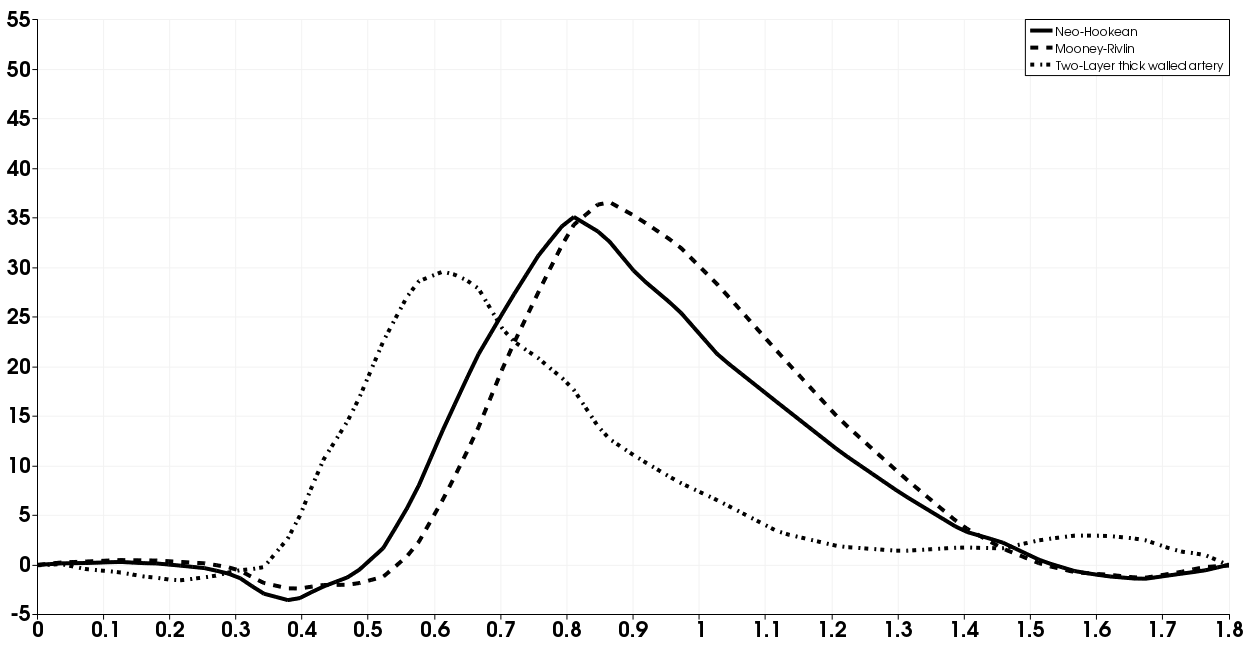}
      \caption{$t=9.00$ ms}
    \end{subfigure}
  }
  \caption{Comparison of fluid pressure waves  along the center line with the 
    starting point $(0, 0, 0)$ cm and ending point $(0, 0, 1.8)$ cm 
    from the FSI simulation using the model of Neo-Hookean material (in solid lines), 
    the Mooney-Rivlin material (in dashed lines) and two-layer thick walled artery (in dash 
    dotted lines).}
  \label{fig:cmp_press_threemodel}
\end{figure}

Secondly, in order to see the effects of different structure models applied in 
the FSI simulation, we also compare the fluid pressure waves extracted from 
the FSI simulation using the model of Neo-Hookean material (solid lines), 
Mooney-Rivlin material (dashed lines) and anisotropic two-layer thick 
walled artery (dash dotted lines) in Fig. \ref{fig:cmp_press_threemodel}, 
where the horizontal line represents the center line (in cm), and the 
vertical line represents the pressure (in Pa). As we observe, 
the simulation results obtained from the model of Neo-Hookean and 
Mooney-Rivlin material are quite similar to each other (the speed and magnitude of 
the pressure waves). This is due to the fact that these two models have 
only one term difference in the energy functional; 
see (\ref{eq:nhenerg}) and (\ref{eq:mrenerg}). The pressure waves obtained 
from the model of the anisotropic two-layer thick walled artery travels 
with slower speed and smaller magnitude than the other two models. 

As we discussed in \cite{ULHY13}, for the partitioned approach, 
we need around $50\sim 55$ fixed-point iterations at each time step; for the 
monolithic approach we need about $4$ Newton iterations. In each fixed-point 
iteration, we need about $4-5$ Newton iterations for solving the fluid 
and structure sub-problems;  and in each Newton iteration, we apply the AMG 
sub-problem solvers for the linearized systems. For each Newton iteration in the 
monolithic approach, we need about $10$ coupled AMG or AMLI iterations; 
and each coupled AMG or AMLI iteration requires apply one iteration of 
AMG sub-problem solvers. Altogether we observe almost $50\%$ saving 
of the computational cost in the monolithic approach in comparison with 
the partitioned approach. Further reduction in computational cost will be 
realized by using parallel computing, see, e.g, \cite{CCDGHUL03}, 
that is considered as a forthcoming work. 

\section{Conclusions}\label{sec:con}
In this work, we have developed the monolithic approach for solving 
the coupled FSI problem in an all-at-once manner. 
The Newton method for the 
nonlinear coupled system demonstrates its robustness and efficiency. For 
solving the linearized FSI system, 
the preconditioned Krylov sub-space, algebraic multigrid and algebraic 
multilevel methods have shown 
their good performance and robustness. In particular, the 
monolithic AMG and AMLI methods show more robustness than 
the preconditioned Krylov sub-space methods utilizing block factorization 
of the coupled system, i.e., the iteration numbers stay in a same range 
with the mesh refinement. 
Compare to the partitioned approach, the monolithic approach developed 
in the work shows its more robustness and efficiency with respect to the 
numerical results and solution methods. 
% latex myarticle
% bibtex myarticle
% latex myarticle
% latex myarticle
\bibliography{FSI_NonLinea}
\bibliographystyle{siam}

\end{document}

%% file: setting.pdf_t
\begin{picture}(0,0)%
\includegraphics{setting.pdf}%
\end{picture}%
\setlength{\unitlength}{3947sp}%
\begingroup\makeatletter\ifx\SetFigFont\undefined%
\gdef\SetFigFont#1#2#3#4#5{%
  \reset@font\fontsize{#1}{#2pt}%
  \fontfamily{#3}\fontseries{#4}\fontshape{#5}%
  \selectfont}%
\fi\endgroup%
\begin{picture}(7980,3019)(1711,-5234)
\put(5851,-2611){\makebox(0,0)[lb]{\smash{{\SetFigFont{12}{14.4}{\rmdefault}{\mddefault}{\updefault}{\color[rgb]{0,0,0}${\mathcal L}^t(\cdot)$}%
}}}}
\put(5851,-4561){\makebox(0,0)[lb]{\smash{{\SetFigFont{12}{14.4}{\rmdefault}{\mddefault}{\updefault}{\color[rgb]{0,0,0}${\mathcal A}^t(\cdot)$}%
}}}}
\put(2851,-4111){\makebox(0,0)[lb]{\smash{{\SetFigFont{12}{14.4}{\rmdefault}{\mddefault}{\updefault}{\color[rgb]{0,0,0}$\Omega_f^0$}%
}}}}
\put(7126,-4186){\makebox(0,0)[lb]{\smash{{\SetFigFont{12}{14.4}{\rmdefault}{\mddefault}{\updefault}{\color[rgb]{0,0,0}$\Omega^t_f$}%
}}}}
\put(3976,-3511){\makebox(0,0)[lb]{\smash{{\SetFigFont{12}{14.4}{\rmdefault}{\mddefault}{\updefault}{\color[rgb]{0,0,0}$\Gamma^0$}%
}}}}
\put(8101,-3586){\makebox(0,0)[lb]{\smash{{\SetFigFont{12}{14.4}{\rmdefault}{\mddefault}{\updefault}{\color[rgb]{0,0,0}$\Gamma^t$}%
}}}}
\put(2026,-3736){\makebox(0,0)[lb]{\smash{{\SetFigFont{12}{14.4}{\rmdefault}{\mddefault}{\updefault}{\color[rgb]{0,0,0}$\Gamma_{in}^0$}%
}}}}
\put(5176,-3736){\makebox(0,0)[lb]{\smash{{\SetFigFont{12}{14.4}{\rmdefault}{\mddefault}{\updefault}{\color[rgb]{0,0,0}$\Gamma_{out}^0$}%
}}}}
\put(9676,-3736){\makebox(0,0)[lb]{\smash{{\SetFigFont{12}{14.4}{\rmdefault}{\mddefault}{\updefault}{\color[rgb]{0,0,0}$\Gamma_{out}^t$}%
}}}}
\put(6451,-3736){\makebox(0,0)[lb]{\smash{{\SetFigFont{12}{14.4}{\rmdefault}{\mddefault}{\updefault}{\color[rgb]{0,0,0}$\Gamma_{in}^t$}%
}}}}
\put(3601,-5161){\makebox(0,0)[lb]{\smash{{\SetFigFont{12}{14.4}{\rmdefault}{\mddefault}{\updefault}{\color[rgb]{0,0,0}$\Omega^0$}%
}}}}
\put(8176,-5161){\makebox(0,0)[lb]{\smash{{\SetFigFont{12}{14.4}{\rmdefault}{\mddefault}{\updefault}{\color[rgb]{0,0,0}$\Omega^t$}%
}}}}
\put(4501,-2461){\makebox(0,0)[lb]{\smash{{\SetFigFont{12}{14.4}{\rmdefault}{\mddefault}{\updefault}{\color[rgb]{0,0,0}$\Omega_s^0$}%
}}}}
\put(9001,-2386){\makebox(0,0)[lb]{\smash{{\SetFigFont{12}{14.4}{\rmdefault}{\mddefault}{\updefault}{\color[rgb]{0,0,0}$\Omega_s^t$}%
}}}}
\put(2776,-2386){\makebox(0,0)[lb]{\smash{{\SetFigFont{12}{14.4}{\rmdefault}{\mddefault}{\updefault}{\color[rgb]{0,0,0}$\Gamma_{n}^0$}%
}}}}
\put(1876,-2536){\makebox(0,0)[lb]{\smash{{\SetFigFont{12}{14.4}{\rmdefault}{\mddefault}{\updefault}{\color[rgb]{0,0,0}$\Gamma_d^0$}%
}}}}
\put(1726,-3511){\makebox(0,0)[lb]{\smash{{\SetFigFont{12}{14.4}{\rmdefault}{\mddefault}{\updefault}{\color[rgb]{0,0,0}$\Gamma_f^0$}%
}}}}
\end{picture}%

%% file: gageo.pdf_t
\begin{picture}(0,0)%
\includegraphics{gageo.pdf}%
\end{picture}%
\setlength{\unitlength}{4144sp}%
\begingroup\makeatletter\ifx\SetFigFont\undefined%
\gdef\SetFigFont#1#2#3#4#5{%
  \reset@font\fontsize{#1}{#2pt}%
  \fontfamily{#3}\fontseries{#4}\fontshape{#5}%
  \selectfont}%
\fi\endgroup%
\begin{picture}(9284,12309)(1982,-11323)
\put(11206,-3211){\makebox(0,0)[lb]{\smash{{\SetFigFont{25}{30.0}{\rmdefault}{\bfdefault}{\updefault}{\color[rgb]{0,0,0}$H_M = 0.26$ mm}%
}}}}
\put(11251,-4111){\makebox(0,0)[lb]{\smash{{\SetFigFont{25}{30.0}{\rmdefault}{\bfdefault}{\updefault}{\color[rgb]{0,0,0}$H_A = 0.13$ mm}%
}}}}
\put(11251,-5011){\makebox(0,0)[lb]{\smash{{\SetFigFont{25}{30.0}{\rmdefault}{\bfdefault}{\updefault}{\color[rgb]{0,0,0}$\alpha_M = 29.0^\circ$ mm}%
}}}}
\put(11251,-6811){\makebox(0,0)[lb]{\smash{{\SetFigFont{25}{30.0}{\rmdefault}{\bfdefault}{\updefault}{\color[rgb]{0,0,0}$R_i = 1.43$ mm}%
}}}}
\put(11251,-7711){\makebox(0,0)[lb]{\smash{{\SetFigFont{25}{30.0}{\rmdefault}{\bfdefault}{\updefault}{\color[rgb]{0,0,0}$L = 18$ mm}%
}}}}
\put(3421,659){\makebox(0,0)[lb]{\smash{{\SetFigFont{25}{30.0}{\rmdefault}{\bfdefault}{\updefault}{\color[rgb]{0,0,0}$R_o$}%
}}}}
\put(4951,-16){\makebox(0,0)[lb]{\smash{{\SetFigFont{25}{30.0}{\rmdefault}{\bfdefault}{\updefault}{\color[rgb]{0,0,0}$R_i$}%
}}}}
\put(3601,-16){\makebox(0,0)[lb]{\smash{{\SetFigFont{25}{30.0}{\rmdefault}{\bfdefault}{\updefault}{\color[rgb]{0,0,0}$H_M$}%
}}}}
\put(2296,-16){\makebox(0,0)[lb]{\smash{{\SetFigFont{25}{30.0}{\rmdefault}{\bfdefault}{\updefault}{\color[rgb]{0,0,0}$H_A$}%
}}}}
\put(2656,-8611){\makebox(0,0)[lb]{\smash{{\SetFigFont{25}{30.0}{\rmdefault}{\bfdefault}{\updefault}{\color[rgb]{0,0,0}$\textup{Adventitia}$}%
}}}}
\put(3556,-5011){\makebox(0,0)[lb]{\smash{{\SetFigFont{25}{30.0}{\rmdefault}{\bfdefault}{\updefault}{\color[rgb]{0,0,0}$\textup{Media}$}%
}}}}
\put(9991,-5191){\makebox(0,0)[lb]{\smash{{\SetFigFont{25}{30.0}{\rmdefault}{\bfdefault}{\updefault}{\color[rgb]{0,0,0}$L$}%
}}}}
\put(7246,-5461){\makebox(0,0)[lb]{\smash{{\SetFigFont{25}{30.0}{\rmdefault}{\bfdefault}{\updefault}{\color[rgb]{0,0,0}$2\alpha_M$}%
}}}}
\put(7246,-9466){\makebox(0,0)[lb]{\smash{{\SetFigFont{25}{30.0}{\rmdefault}{\bfdefault}{\updefault}{\color[rgb]{0,0,0}$2\alpha_A$}%
}}}}
\put(11251,-5911){\makebox(0,0)[lb]{\smash{{\SetFigFont{25}{30.0}{\rmdefault}{\bfdefault}{\updefault}{\color[rgb]{0,0,0}$\alpha_A = 62.0^\circ$ mm}%
}}}}
\end{picture}%